\documentclass[twocolumn,twoside]{IEEEtran}
\usepackage{graphicx}
\usepackage{amsmath,amssymb,amsfonts,color,psfrag}
\usepackage{algorithm,algorithmic}

\newtheorem{theorem}{Theorem}[section]
\newtheorem{lemma}[theorem]{Lemma}
\newtheorem{remark}[theorem]{Remark}
\newtheorem{corollary}[theorem]{Corollary}
\newtheorem{proposition}[theorem]{Proposition}

\newcommand{\real}{{\mathbb{R}}}

\renewcommand{\natural}{{\mathbb{N}}}
\newcommand{\eps}{\epsilon}

\newcommand{\pder}[2]{\frac{\partial #1}{\partial #2}}

\newcommand{\dist}{\operatorname{dist}}
\newcommand{\HH}{{\mathcal{H}}}
\newcommand{\VV}{{\mathcal{V}}}
\newcommand{\WW}{{\mathcal{W}}}
\newcommand{\NN}{{\mathcal{N}}}

\newcommand{\SR}{\operatorname{SR}}

\newcommand{\subscr}[2]{#1_{\text{#2}}}
\newcommand{\setdef}[2]{\{#1 \; | \enspace #2\}}

\newcommand{\ov}{\overline}

\newcommand{\algoASR}{\textsc{Adjust sensing radius algorithm}}
\newcommand{\algoACR}{\textsc{Adjust communication radius algorithm}}
\newcommand{\algoM}{\textsc{Monitoring algorithm}}
\newcommand{\algoCone}{\textsc{Coverage behavior algorithm I}}
\newcommand{\algoCtwo}{\textsc{Coverage behavior algorithm II}}

\renewcommand{\theenumi}{(\roman{enumi})}

\begin{document}
\title{Coverage control for mobile sensing networks\thanks{Submitted on
    November 4, 2002.  Previous short versions of this paper appeared in the
    IEEE Conference on Robotics and Automation, Arlington, VA, May 2002, and
    Mediterranean Conference on Control and Automation, Lisbon, Portugal, July
    2002.}}

\author{Jorge Cort\'es, \thanks{Jorge Cort\'es, Timur Karatas and Francesco
    Bullo are with the Coordinated Science Laboratory, University of Illinois
    at Urbana-Champaign, 1308 W.  Main St., Urbana, IL 61801, United States,
    Tels: +1-217-244-8734, +1-217-244-9414 and +1-217-333-0656, Fax:
    +1-217-244-1653, Email: \texttt{\{jcortes,tkaratas,bullo\}@uiuc.edu}}
  Sonia Mart{\'\i}nez, \thanks{Sonia Mart{\'\i}nez is with the Escola
    Universit\`aria Polit\`ecnica de Vilanova i la Geltr\'u, Universidad
    Polit\'ecnica de Catalu\~na, Av.  V.  Balaguer s/n, Vilanova i la
    Geltr\'u, 08800, Spain, Tel: +34-938967720, Fax: +34-938967700, Email:
    \texttt{soniam@mat.upc.es}} Timur Karatas, Francesco Bullo, {\it Member
    IEEE}}

\markboth{Submitted as a regular paper to IEEE Transactions on Robotics
  and Automation}{Cort\'es, Mart{\'\i}nez, Karatas and Bullo: Coverage
  control for mobile sensing networks}

\maketitle

\begin{abstract}
  This paper presents control and coordination algorithms for groups of
  vehicles.  The focus is on autonomous vehicle networks performing
  distributed sensing tasks where each vehicle plays the role of a mobile
  tunable sensor.  The paper proposes gradient descent algorithms for a class
  of utility functions which encode optimal coverage and sensing policies.
  The resulting closed-loop behavior is adaptive, distributed, asynchronous,
  and verifiably correct.
\end{abstract}


\section{Introduction}

\subsection*{Mobile sensing networks}

The deployment of large groups of autonomous vehicles is rapidly becoming
possible because of technological advances in networking and in
miniaturization of electro-mechanical systems.  In the near future large
numbers of robots will coordinate their actions through ad-hoc communication
networks and will perform challenging tasks including search and recovery
operations, manipulation in hazardous environments, exploration,
surveillance, and environmental monitoring for pollution detection and
estimation.  The potential advantages of employing teams of agents are
numerous.  For instance, certain tasks are difficult, if not impossible, when
performed by a single vehicle agent.  Further, a group of vehicles inherently
provides robustness to failures of single agents or communication links.

Working prototypes of active sensing networks have already been developed;
see~\cite{CRW-al:99,EK-JB:99,PER-al:00}.  In~\cite{PER-al:00}, launchable
miniature mobile robots communicate through a wireless network.  The vehicles
are equipped with sensors for vibrations, acoustic, magnetic, and IR signals
as well as an active video module (i.e., the camera or micro-radar is
controlled via a pan-tilt unit).  A second system is suggested
in~\cite{TBC-JGB-JC-DW:93} under the name of Autonomous Oceanographic
Sampling Network; see also~\cite{RMT-EHT:98,EE-SP:99}.  In this case,
underwater vehicles are envisioned measuring temperature, currents, and other
distributed oceanographic signals.  The vehicles communicate via an acoustic
local area network and coordinate their motion in response to local sensing
information and to evolving global data.  This mobile sensing network is
meant to provide the ability to sample the environment adaptively in space
and time.  By identifying evolving temperature and current gradients with
higher accuracy and resolution than current static sensors, this technology
could lead to the development and validation of improved oceanographic
models.

\subsection*{Optimal sensor allocation and coverage problems}

A fundamental prototype problem in this paper is that of characterizing and
optimizing notions of quality-of-service provided by an adaptive sensor
network in a dynamic environment.  To this goal, we introduce a notion of
\emph{sensor coverage} that formalizes an optimal sensor placement problem.
This spatial resource allocation problem is the subject of a discipline
called locational
optimization~\cite{AO-BB-KS:94,ZD:95,AS-AO:95,AO-AS:97,AO-BB-KS-SNC:00}.

Locational optimization problems pervade a broad spectrum of scientific
disciplines. Biologists rely on locational optimization tools to study how
animals share territory and to characterize the behavior of animal groups
obeying the following interaction rule: each animal establishes a region of
dominance and moves toward its center. Locational optimization problems are
spatial resource allocation problems (where to place mailboxes in a city or
cache servers on the internet) and play a central role in quantization and
information theory (the design of a minimum-distortion fixed-rate vector
quantizer is a locational problem).  Other technologies affected by
locational optimization include mesh and grid optimization methods,
clustering analysis, data compression, and statistical pattern recognition.

Because locational optimization problems are so widely studied, it is not
surprising that methods are indeed available to tackle coverage problems;
see~\cite{AO-BB-KS:94,AO-AS:97,QD-VF-MG:99,AO-BB-KS-SNC:00}.  However, most
currently-available algorithms are not applicable to mobile sensing networks
because they inherently assume a centralized computation for a limited size
problem in a known static environment.  This is not the case in multi-vehicle
networks which, instead, rely on a distributed communication and computation
architecture. Although an ad-hoc wireless network provides the ability to
share some information, no global omniscient leader might be present to
coordinate the group.  The inherent spatially-distributed nature and limited
communication capabilities of a mobile network invalidate classic approaches
to algorithm design.

%% For example, the vehicle ensemble is endowed
% with a communication network with changing topology, limited bandwidth,
% delays, and faulty links.  In the physical domain, we consider networks
% dynamically varying in size and composition.  It might be the case that the
% vehicles' motion needs to respect certain relative distances or line of sight
% constraints.  Furthermore, we design algorithms that provide guaranteed
% quality of service.

\subsection*{Distributed asynchronous algorithms for coverage control}

In this paper we design coordination algorithms implementable by a
multi-vehicle network with limited sensing and communication capabilities.
Our approach is related to the classic Lloyd algorithm from quantization
theory; see~\cite{SPL:82} for a reprint of the original report and
\cite{RMG-DLN:98} for a historical overview.  We present Lloyd descent
algorithms that take into careful consideration all constraints on the mobile
sensing network.  In particular, we design coverage algorithms that are
adaptive, distributed, asynchronous, and verifiably asymptotically correct:
\begin{description}
\item[Adaptive:] Our coverage algorithms provide the network with
  the ability to address changing environments, sensing task, and network
  topology (due to agents departures, arrivals, or failures).
\item[Distributed:] Our coverage algorithms are distributed in the sense that
  the behavior of each vehicle depends only on the location of its neighbors.
  Also, our algorithms do not required a fixed-topology communication graph,
  i.e., the neighborhood relationships do change as the network evolves.  The
  advantages of distributed algorithms are scalability and robustness.
\item[Asynchronous:] Our coverage algorithms are amenable to asynchronous
  implementation.  This means that the algorithms can be implemented in a
  network composed of agents evolving at different speeds, with different
  computation and communication capabilities.  Furthermore, our algorithms do
  not require a global synchronization and convergence properties are
  preserved even if information about neighboring vehicles propagates with
  some delay.  An advantage of asynchronism is a minimized communication
  overhead.
\item[Verifiable Asymptotically Correct:] Our algorithms guarantees monotonic
  descent of the cost function encoding the sensing task. Asymptotically the
  evolution of the mobile sensing network is guaranteed to converge to
  so-called centroidal Voronoi configurations that are critical points of the
  optimal sensor coverage problem.  
%  The importance of formal verification proofs will increase with the
%  dimension and complexity of vehicle networks.
\end{description}

Let us describe in some detail what are the contributions of this paper.
Section~\ref{sec:review} reviews certain locational optimization problems and
their solutions as centroidal Voronoi partitions.
Section~\ref{sec:coverage-control} provides a continuous-time version of the
classic Lloyd algorithm from vector quantization and applies it to the
setting of multi-vehicle networks. In discrete-time, we propose a family of
Lloyd algorithms.  We carefully characterize convergence properties for both
continuous and discrete-time versions (Appendix~\ref{sec:appendix} collects
some relevant facts on descent flows).  We discuss a worst-case optimization
problem, we investigate a simple uniform planar setting, and we present
numerical results.

Section~\ref{sec:distributed-asynchronous} presents two asynchronous
distributed implementations of Lloyd algorithm for ad-hoc networks with
communication and sensing capabilities. Our treatment carefully accounts for
the constraints imposed by the distributed nature of the vehicle network.  We
present two asynchronous implementations, one based on classic results on
distributed gradient flows, the other based on the structure of the coverage
problem.

Section~\ref{sec:vehicle-dynamics} considers vehicle models with more
realistic dynamics.  We present two formal results on passive vehicle
dynamics and on vehicles equipped with individual local controllers. We
present numerical simulations of passive vehicle models and of unicycle
mobile vehicles. Next, Section~\ref{sec:geometric-patterns} describes density
functions that lead the multi-vehicle network to predetermined geometric
patterns.  

%Finally, we present our conclusions and directions for future research in
%Section~\ref{sec:conclusions}.

\subsection*{Review of distributed algorithms for cooperative control}
Recent years have witnessed a large research effort focused on motion
planning and coordination problems for multi-vehicle systems.  Issues include
geometric patterns \cite{HY-TA:94,KS-IS:96}, formation
control~\cite{TB-RA:98,ME-XH:01b,JPD-JPO-VK:01,PT-GP-PL:02,ROS-RMM:02b},
%% cooperative motion planning~\cite{ME-TLP:87,SML-SAH:98b}, 
and conflict avoidance~\cite{CT-GJP-SSS:98,EF-ZHM-JHO-EF:01}.  Algorithms for
robotic sensing tasks are presented for example in~\cite{HC:01,RB-NEL:02}.
It is only recently, however, that truly distributed coordination laws for
dynamic networks are being proposed; e.g., see~\cite{AJ-JL-ASM:02,EK:02a} and
the conference versions of this work~\cite{JC-SM-TK-FB:01k,JC-SM-TK-FB:02b}.

Heuristic approaches to the design of interaction rules and emerging
behaviors have been throughly investigated within the literature on
behavior-based robotics;
see~\cite{RAB:86,CWR:87,RCA:98,TB-RA:98,MSF-MJM:98,ACS-LEP:02,TB-LEP:02,LEP:02a}.
An example of coverage control is discussed in~\cite{AH-MJM-GSS:02}. Along
this line of research, algorithms have been designed for sophisticated
cooperative tasks. However, no formal results are currently available on how
to design reactive control laws, ensure their correctness, and guarantee
their optimality with respect to an aggregate objective.

The study of distributed algorithms is concerned with providing mathematical
models, devising precise specifications for their behavior, and formally
proving their correctness and complexity.  Via an automata-theoretic
approach, the references~\cite{NAL:97,GT:01} treat distributed consensus,
resource allocation, communication, and data consistency problems.  From a
numerical optimization viewpoint, the works
in~\cite{JNT-DPB-MA:86,DPB-JNT:97,SHL-DEL:99} discuss distributed
asynchronous algorithms as networking algorithms, rate and flow control, and
gradient descent flows. Typically, both these sets of references consider
networks with fixed topology, and do not address algorithms over ad-hoc
dynamically changing networks. Another common assumption is that, any time an
agent communicates its location, it broadcasts it to every other agent in the
network. In our setting, this would require a non-distributed communication
set-up.

\section{From location optimization to centroidal Voronoi partitions}
\label{sec:review}
\subsection{Locational optimization}
In this section we describe a collection of known facts about a meaningful
optimization problem.  References include the theory and applications of
centroidal Voronoi partitions, see~\cite{QD-VF-MG:99}, and the discipline of
facility location, see~\cite{ZD:95}. Along the paper, we interchangeably
refer to the elements of the network as sensors, agents, vehicles, or robots.

Let $Q$ be a convex polytope in $\real^N$ and let $\|\cdot\|$\ denote the
Euclidean distance function.  We call a map $\phi:Q\rightarrow\real_+$\ a
\emph{distribution density function} if it represents a measure of
information or probability that some event take place over $Q$. In equivalent
words, we can consider $Q$ to be the bounded support of the function $\phi$.
Let $P=(p_1,\dots,p_n)$ be the \emph{location of~$n$ sensors}, each moving in
the space~$Q$.  Because of noise and loss of resolution, the {\emph{sensing
    performance}} at point~$q$ taken from $i$th sensor at the position $p_i$
degrades with the distance $\|q -p_i \|$ between $q$ and $p_i$; we describe
this degradation with a non-decreasing differentiable function $f:\real_+ \to
\real_+$.  Accordingly, $f\left(\|q - p_i \| \right)$ provides a quantitative
assessment of how poor the sensing performance is.

\begin{remark}
  As an example, consider $n$ mobile robots equipped with microphones
  attempting to detect, identify, and localize a sound-source.  \emph{How
    should we plan to robots' motion in order to maximize the detection
    probability?}  Assuming the source emits a known signal, the optimal
  detection algorithm is a matched filter (i.e., convolve the known waveform
  with the received signal and threshold).  The source is detected depending
  on the signal-to-noise-ratio, which is inversely proportional to the
  distance between the microphone and the source.  Various electromagnetic
  and sound sensors have signal-to-noise ratios inversely proportional to
  distance.
\end{remark}

Within the context of this paper, a \emph{partition} of $Q$
is a collection of $n$ polytopes $\WW=\{W_1,\dots,W_n\}$
with disjoint interiors whose union is $Q$.  We say that two
partitions $\WW$ and $\WW'$ are equal if $W_i$ and $W_i'$
only differ by a set of $\phi$-measure zero, for all $i \in
\{ 1,\dots,n\}$.

We consider the task of minimizing the locational
optimization function
\begin{equation} \label{eq:HH}
  \HH (P,\WW) = \sum_{i=1}^n \int_{W_i} f(\| q-p_i \|) d\phi(q),  
\end{equation}
where we assume that the $i$th sensor is responsible for
measurements over its ``dominance region''~$W_i$.  Note that
the function $\HH$ is to be minimized with respect to both
(1) the sensors location $P$, and (2) the assignment of the
dominance regions $\WW$. This problem is referred to as a
facility location problem and in particular as a continuous
$p$-median problem in~\cite{ZD:95}.

\begin{remark}
  Note that if we interchange the positions of any two
  agents, along with their associated regions of dominance,
  the value of the locational optimization function $\HH$ is
  not affected.  To eliminate this discrete redundancy, one
  could take the discrete group of permutations $\Sigma_n$
  with the natural action on $Q^n$, and consider $Q^n /
  \Sigma_n$ as the configuration space for the position $P$
  of the $n$ vehicles.
\end{remark}

\subsection{Voronoi partitions}
One can easily see that, at fixed sensors location, the
optimal partition of $Q$ is the \emph{Voronoi partition}
$\VV(P)=\{V_1,\dots,V_n\}$ generated by the points
$(p_1,\dots,p_n)$:
\begin{equation*}
  V_i = \setdef{q\in Q}{\| q - p_i\| \leq \|q - p_j \| \, , \; \forall j\neq i}.
\end{equation*}
We refer to~\cite{AO-BB-KS-SNC:00} for a comprehensive treatment on Voronoi
diagrams, and briefly present some relevant concepts.  The set of regions
$\{V_1,\dots,V_n\}$\ is called the Voronoi diagram for the generators
$\{p_1,\dots,p_n\}$.  When the two Voronoi regions $V_i$\ and $V_j$\ are
adjacent, $p_i$ is called a \emph{(Voronoi) neighbor} of $p_j$ (and
vice-versa).  The set of indexes of the Voronoi neighbors of $p_i$ is denoted
by $\NN (i)$.  Clearly, $j \in \NN (i)$ if and only if $i \in \NN (j)$.  We
also define the $(i,j)$-face as $\Delta_{ij}=V_i\cap V_j$.  Voronoi diagrams
can be defined with respect to various distance functions, e.g., the $1$-,
$2$-, $s$-, and $\infty$-norm over $Q=\real^m$, and Voronoi diagrams can be
defined over Riemannian manifolds; see~\cite{GL-DL:00}.  Some useful facts
about the Euclidean setting are the following: if $Q$ is a convex polytope in
a $N$-dimensional Euclidean space, the boundary of each $V_i$\ is the union
of $(N-1)$-dimensional convex polytopes.

In what follows, we shall write
\begin{equation*}
  \HH_{\VV} (P)= \HH(P,\VV(P)).
\end{equation*}
Note that
\begin{align}\label{eq:HH_VV}
  \HH_{\VV} (P) &= \int_Q \min_{i\in\{1,\ldots,n\}} f(\| q - p_i \| ) d\phi (q) \, ,
  \\
  &= E\left[  \min_{i\in\{1,\ldots,n\}} f(\|q-p_i\|) \right],
  \nonumber
\end{align}
that is, the locational optimization function can be interpreted as an
expected value composed with a min operation. This is the usual way in which
the problem is presented in the facility location and operations research
literature~\cite{ZD:95,AS-AO:95}.  Remarkably, one can
show~\cite{QD-VF-MG:99} that
\begin{align} \label{eq:remarkable-derivative}
  \pder{\HH_{\VV}}{p_i}(P) &= \pder{\HH}{p_i}(P,\VV(P))
  = \int_{V_i} \pder{}{p_i} f\left(\| q - p_i \| \right) d\phi(q),
\end{align}
and deduce some smoothness properties of $\HH_{\VV}$.  Since
the Voronoi partition $\VV$ depends at least continuously on
$P=(p_1,\ldots,p_n)$, the function $\HH_{\VV}$ is at least
continuously differentiable.  

\subsection{Centroidal Voronoi partitions}
Let us recall some basic quantities associated to a region
$V\subset\real^N$ and a mass density function $\rho$.  The
(generalized) mass, centroid (or center of mass), and polar
moment of inertia are defined as
\begin{gather*} 
  M_V = \int_V \rho(q) \, dq, \quad C_V = \frac{1}{M_V}
  \int_{V} q \, \rho(q) \, dq,
  \\
  J_{V,p} = \int_{V} \|q-p\|^2 \, \rho(q) \, dq.
\end{gather*}
Additionally, by the parallel axis theorem, one can write,
\begin{align}\label{eqn:jpc}
  J_{V,p} &= J_{V,C_V} + M_V \, \|p-C_V\|^2
\end{align}
where $J_{V,C_V}\in\real_+$ is defined as the polar moment
of inertia of the region $V$ about its centroid.

Let us consider again the locational optimization
problem~\eqref{eq:HH}, and suppose now we are strictly
interested in the setting
\begin{equation} \label{eq:HH-norm}
  \HH (P,\WW) = \sum_{i=1}^n \int_{W_i} \|q-p_i\|^2 d\phi(q),  
\end{equation}
that is, we assume $f(\|q-p_i\|)=\|q-p_i\|^2$.  Applying the
parallel axis theorem leads to simplifications for both the
function $\HH_{\VV}$ and its partial derivative:
\begin{align*}
  \HH_{\VV}(P)  &=
  \sum_{i=1}^n J_{V_i,C_{V_i}}+
  \sum_{i=1}^n M_{V_i} \,  \|p_i-C_{V_i}\|^2  \\
  \pder{\HH_{\VV}}{p_i}(P) &=  2 M_{V_i} (p_i - C_{V_i}) .
\end{align*}
It is convenient to define $\HH_{\VV,1} =\sum_{i=1}^n
J_{V_i,C_{V_i}}$ and $\HH_{\VV,2} =\sum_{i=1}^n M_{V_i}
\|p_i-C_{V_i}\|^2$.

Therefore, the (not necessarily unique) local minimum points
for the location optimization function $\HH_{\VV}$ are
{\emph{centroids}} of their Voronoi cells, i.e., each
location $p_i$ satisfies two properties simultaneously: it
is the generator for the Voronoi cell $V_i$ and it is its
centroid
\begin{align*}
  C_{V_i} & = \operatorname{argmin}_{p_i} \HH_{\VV}(P).
\end{align*}
Accordingly, the critical partitions and points for $\HH$
are called \emph{centroidal Voronoi partitions}. We will
refer to a sensors configuration as a \emph{centroidal
  Voronoi configuration} if it gives rise to a centroidal
Voronoi partition.  This discussion provides a proof
alternative to the one given in~\cite{QD-VF-MG:99} for the
necessity of centroidal Voronoi partitions as solutions to
the continuous $p$-median location problem.

\section{Continuous and Discrete-Time Lloyd Descent for
  Coverage Control} 
\label{sec:coverage-control}
In this section, we describe algorithms to compute the location of sensors
that minimize the cost $\HH$, both in continuous and in discrete-time.  In
Section~\ref{se:continuous-Lloyd}, we propose a continuous-time version of
the classic Lloyd algorithm. Here, both the positions and partitions evolve
in continuous time, whereas Lloyd algorithm for vector quantization is
designed in discrete time.  In Section~\ref{se:discrete-Lloyd}, we develop a
family of variations of Lloyd algorithm in discrete time.  In both setting,
we prove that the proposed algorithms are \emph{gradient descent flows}.

\subsection{A continuous-time Lloyd algorithm}\label{se:continuous-Lloyd}
Assume the sensors location obeys a first order dynamical
behavior described by
\begin{equation*}
  \dot{p}_i = u_i.
\end{equation*}
Consider $\HH_{\VV}$ a cost function to be minimized and
impose that the location $p_i$ follows a gradient descent.
In equivalent control theoretical terms, consider
$\HH_{\VV}$ a Lyapunov function and stabilize the
multi-vehicle system to one of its local minima via
dissipative control.  Formally, we set
\begin{gather} \label{eq:continuous-lloyd}
  u_i=  - \subscr{k}{prop} (p_i - C_{V_i}),
\end{gather}
where $k$ is a positive gain, and where we assume that the
partition $\VV(P)=\{V_1,\dots,V_n\}$ is continuously
updated.
\begin{proposition}[Continuous-time Lloyd descent]
  For the closed-loop system induced by equation~\eqref{eq:continuous-lloyd},
  the sensors location converges asymptotically to the set of critical points
  of~$\HH_\VV$, i.e., the set of centroidal Voronoi configurations on $Q$.
  Assuming this set is finite, the sensors location converges to a centroidal
  Voronoi configuration.
\end{proposition}
\begin{proof}
  Under the control law~\eqref{eq:continuous-lloyd}, we have
  \begin{multline*}
    \frac{d}{dt} \HH_{\VV} (P(t)) = \sum_{i=1}^n \pder{\HH_{\VV}}{p_i}\,
    \dot{p}_i \\
    = -2 \subscr{k}{prop} \sum_{i=1}^n M_{V_i} \| p_i - C_{V_i} \|^2 = - 2
    \subscr{k}{prop} \HH_{\VV,2} (P(t)).
  \end{multline*}
  By LaSalle's principle, the sensors location converges to
  the largest invariant set contained in
  $\HH_{\VV,2}^{-1}(0)$, which is precisely the set of
  centroidal Voronoi configurations.  Since this set is
  clearly invariant for~\eqref{eq:continuous-lloyd}, we get
  the stated result.  If $\HH_{\VV,2}^{-1}(0)$ consists of a
  finite collection of points, then $P(t)$ converges to one
  of them, see Corollary~\ref{corollary:LaSalle-finite}.
\end{proof}

\begin{remark}
  If $\HH_{\VV,2}^{-1}(0)$ is finite, and $P(t) \rightarrow C$, then a
  sufficient condition that guarantees exponential convergence is that the
  Hessian of $\HH_\VV$ be positive definite at $C$. This property is known to
  be an open problem, see~\cite{QD-VF-MG:99}.  Note that this gradient
  descent is not guaranteed to find the global minimum.  For example, in the
  vector quantization and signal processing literature~\cite{RMG-DLN:98}, it
  is known that for bimodal distribution density functions, the solution to
  the gradient flow reaches local minima where the number of generators
  allocated to the two region of maxima are not optimally partitioned.
\end{remark}

\subsection{A family of discrete-time Lloyd
  algorithms}\label{se:discrete-Lloyd}

Let us consider the following class of variations of Lloyd
algorithm.  Let $T$ be a continuous mapping $T : Q^n
\rightarrow Q^n$ verifying the following two properties:
\begin{description}
\item (a) for all $i \in \{1, \dots, n \}$, $\| T_i(P) -
  C_{V_i(P)} \| \le \| p_i - C_{V_i(P)} \|$, where $T_i$
  denotes the $i$th component of $T$,
\item (b) if $P$ is not centroidal, then there exists a $j$
  such that $\| T_j(P) - C_{V_j(P)} \| < \| p_j - C_{V_j(P)}
  \|$.
\end{description}
Property (a) guarantees that, if moving, the agents of the
network do not increase their distance to its corresponding
centroid.  Property (b) ensures that at least one robot
moves at each iteration and strictly approaches the centroid
of its Voronoi region. Because of this property, the fixed
points of $T$ are the set of centroidal Voronoi
configurations.

% Typically one could also add a condition of the form:
% \begin{description}
% \item (c)  there exists  a fixed  $m \in \natural$  such that,  if $P$  is
% not centroidal, then $T^m_i(P) \not = p_i$ for all $p_i \not = C_{V_i}$. 
% \end{description}
% That is, after a  fixed number of iterations, all agents which  are not at
% the centroid position of its Voronoi cell have moved. However, property (c)
% is not necessary to prove the convergence of the algorithm $T$, although it
% certainly contributes to speed it up.

\begin{proposition}[Discrete-time Lloyd descent] \label{prop:discrete-Lloyd}
  Let $P_0$ $\in Q^n$ denote the initial sensors location.  Then, the
  sequence $\{T^m(P_0) \}_{m \ge 1}$ converges to the set of centroidal
  Voronoi configurations.  If this set if finite, then $\{T^m(P_0) \}_{m \ge
    1}$ converges to a centroidal Voronoi configuration.
\end{proposition}

\begin{proof}
  Consider $\HH_\VV: Q^n \rightarrow \real_+$ as an
  objective function for the algorithm $T$. Note that
  \begin{gather} \label{eq:inequality-I}
    \HH (P,\VV (P)) \le \HH (P,\WW) \, ,
  \end{gather}
  with strict inequality if $\WW \not = \VV (P)$.  Moreover,
  the parallel axis theorem guarantees
  \begin{gather} \label{eq:inequality-II}
    \HH (P',\WW) \le \HH(P,\WW) \, ,
  \end{gather}
  as long as $\| p_i' - C_{W_i} \| \le \| p_i - C_{W_i} \|$
  for all $i \in \{1,\dots,n \}$, with strict inequality if
  for any $i$, $\| p_i' - C_{W_i} \| < \| p_i - C_{W_i} \|$.
  In particular, $\HH (C_\WW,\WW) \le \HH (P,\WW)$, with
  strict inequality if $P \not = C_\WW$, where $C_\WW$
  denotes the set of centroids of the partition $\WW$.
  
  Now, we have
  \[
  \HH_{\VV} (T(P)) = \HH (T(P), \VV (T(P))) \le \HH (T(P), \VV
  (P)) \, ,
  \]
  because of~\eqref{eq:inequality-I}.  In addition, because
  of property (a) of $T$,
  inequality~\eqref{eq:inequality-II} yields
  \[
  \HH (T(P), \VV (P)) \le \HH (P, \VV (P)) = \HH_{\VV} (P) \,
  ,
  \]
  and the inequality is strict if $P$ is not centroidal by property (b) of
  $T$. Hence, $\HH_\VV$ is a descent function for the algorithm $T$.  The
  result now follows from the global convergence
  Theorem~\ref{lemma:discrete-LaSalle} and
  Proposition~\ref{prop:discrete-LaSalle-surprising}.
\end{proof}

\begin{remark}
  Lloyd algorithm in quantization theory~\cite{SPL:82,RMG-DLN:98} is usually
  presented as follows: given the location of $n$ agents, $p_1, \dots, p_n$,
  (i) construct the Voronoi partition corresponding to $P=(p_1,\dots,p_n)$;
  (ii) compute the mass centroids of the Voronoi regions found in step (i).
  Set the new location of the agents to these centroids; and return to step
  (i).  Lloyd algorithm can also be seen as a fixed point iteration.
  Consider the mappings $LL_i : Q^n \rightarrow Q$ for ${i\in\{1,\ldots,n\}}$
  \[
  LL_i(p_1,\dots,p_n) = \left({\int_{V_i(P)} \phi(q) dq}\right)^{-1}
  {\int_{V_i(P)} q \phi(q) dq} \, .
  \]
  Let $LL: Q^n \rightarrow Q^n$ be defined by $LL = (LL_1, \dots,LL_n)$.
  Clearly, $LL$ is continuous (indeed, $C^1$), and corresponds to Lloyd
  algorithm.  Now, $\| LL_i (P) - C_{V_i} \| = 0 \le \| p_i - C_{V_i} \|$,
  for all $i \in \{ 1, \dots, n \}$.  Moreover, if $P$ is not centroidal,
  then the inequality is strict for all $p_i \not = C_{V_i}$.  Therefore,
  $LL$ verifies properties (a) and~(b).
  
\end{remark}

\subsection{Generalized settings, worst-case design, and the $p$-center
  problem}

Different sensor performance functions $f$ in equation~\eqref{eq:HH}
correspond to different optimization problems.  Provided one uses the
Euclidean distance in the definition of $\HH_{\WW}$, the standard Voronoi
partition computed with respect to the Euclidean metric remains the optimal
partition.  For arbitrary $f$, it is not possible anymore to decompose
$\HH_\VV$ into the sum of terms similar to $\HH_{\VV,1}$ and $\HH_{\VV,2}$.
Nevertheless, it is still possible to implement the gradient flow via the
expression for the partial derivative~\eqref{eq:remarkable-derivative}.
\begin{proposition}
  Assume the sensors location obeys a first order dynamical
  behavior, $\dot{p}_i = u_i$.  Then, for the closed-loop
  system induced by the gradient
  law~\eqref{eq:remarkable-derivative}, $u_i = - {\partial
    \HH_\VV}/{\partial p_i}$, the sensors location
  $P=(p_1,\dots,p_n)$ converges asymptotically to the set of
  critical points of $\HH_\VV$.  Assuming this set is
  finite, the sensors location converges to a critical
  point.
\end{proposition}

More generally, various distance notions can be used to define locational
optimization functions.  Different performance function gives rise to
corresponding notions of ``center of a region'' (any notion of geometric
center, mean, or average is an interesting candidate).  These can then be
adopted in designing coverage algorithms.  We refer to~\cite{RK:89} for a
discussion on Voronoi partitions based on non-Euclidean distance functions
and to~\cite{AO-BB-KS:94,AO-AS:97} for a discussion on the corresponding
locational optimization problems.

Next, let us discuss an interesting variation of the original problem.
In~\cite{ZD:95,AS-AO:95}, minimizing the expected minimum distance function
$\HH_{\VV}$ in equation~\eqref{eq:HH_VV} is referred to as the
\emph{continuous $p$-median problem}.  It is instructive to consider the
worst-case minimum distance function, corresponding to the scenario where no
information is available on the distribution density function.  In other
words, the network seeks to minimize the largest possible distance from any
point in $Q$ to any of the sensor locations, i.e., to minimize the function
\begin{align*}
   \max_{q\in Q} \left[
    \min_{i\in\{1,\ldots,n\}} \|q-p_i\| \right] =
  \max_{i\in\{1,\ldots,n\}} \left[ \max_{q\in V_i}
    \|q-p_i\| \right] \, . 
\end{align*}
This optimization is referred to as the \emph{$p$-center problem}
in~\cite{AS-AO:95,AS-ZD:96}.  One can design a strategy for the $p$-center
problem analog to the Lloyd algorithm for the $p$-median problem: each
vehicle moves, in continuous or discrete-time, toward the center of the
minimum-radius sphere enclosing the polytope.
% Note that the center of the minimum spanning
% sphere for a given convex polytope can be computed via a
% convex problem~\cite{SB-LV:02}.
To the best of our knowledge, no convergence proof is
available in the literature for this algorithm; e.g.,
see~\cite{AS-ZD:96}.  We refer to~\cite{JC-FB:02m} for a
convergence analysis of the continuous and discrete time
algorithms.

In what follows, we shall restrict our attention to the $p$-median problem and
to centroidal Voronoi partitions.

\subsection{Computations over polygons with uniform
  density}

\begin{figure}[htbp]
  \centering
  \includegraphics[width=0.6\linewidth,height=0.4\linewidth]{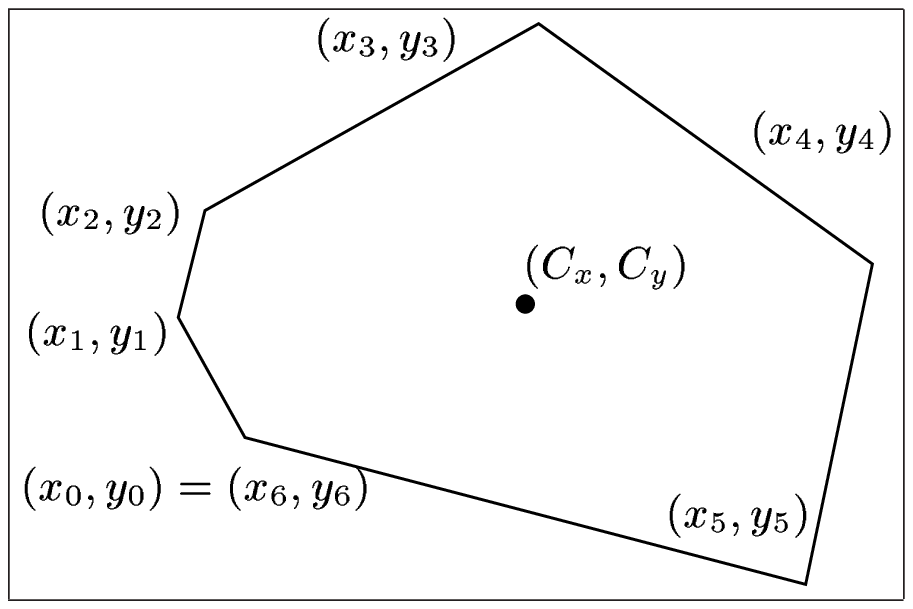} 
  \caption{Notation conventions for a convex polygon.}
  \label{fig:polygon}
\end{figure}

In this section, we investigate closed-form expression for the control laws
introduced above.  Assume the Voronoi region $V_i$ is a convex polygon (i.e.,
a polytope in $\real^2$) with $N_i$ vertexes labeled
$\{(x_0,y_0),\dots,(x_{N_i-1},y_{N_i-1})\}$ such as in
Figure~\ref{fig:polygon}.  It is convenient to define
$(x_{N_i},y_{N_i})=(x_0,y_0)$.  Furthermore, we assume that the density
function is $\phi(q)=1$.
By evaluating the corresponding integrals, one can obtain
the following closed-form expressions
\begin{align}\label{eq:closed-expressions}
  M_{V_i} &= \frac{1}{2}\sum_{k=0}^{N_i-1} (x_k y_{k+1}-x_{k+1}y_k) \nonumber
  \\
  C_{{V_i},x} &= \frac{1}{6M_{V_i}} \sum_{k=0}^{N_i-1} 
  (x_k+x_{k+1})(x_ky_{k+1}-x_{k+1}y_k) \\
  C_{{V_i},y} &= \frac{1}{6M_{V_i}} \sum_{k=0}^{N_i-1}
  (y_k+y_{k+1})(x_ky_{k+1}-x_{k+1}y_k) \, . \nonumber
\end{align}
To present a simple formula for the polar moment of inertia,
let $\bar{x}_k=x_k-C_{{V_i},x}$ and
$\bar{y}_k=y_k-C_{{V_i},y}$, for $k\in\{0,\dots,N_i-1\}$.
Then, the polar moment of inertia of a polygon about its
centroid, $J_{{V_i},C}$ becomes
\begin{multline*}
  J_{{V_i},C_{V_i}} = \frac{1}{12} \sum_{k=0}^{N_i-1}
  (\bar{x}_k\bar{y}_{k+1} - \bar{x}_{k+1}\bar{y}_k) \; \cdot \\
  ( \bar{x}_k^2 + \bar{x}_kx_{k+1} + \bar{x}_{k+1}^2 + \bar{y}_k^2 +
  \bar{y}_k \bar{y}_{k+1} + \bar{y}_{k+1}^2) \,.
\end{multline*}
The proof of these formulas is based on decomposing the polygon into the
union of disjoint triangles.  We refer to~\cite{CC-AP:90} for analog
expressions over~$\real^N$.

A second observation is that the Voronoi polygon's vertexes can be expressed
as a function of the neighboring vehicles.  The vertexes of the $i$th Voronoi
polygon which lie in the interior of $Q$ are the circumcenters of the
triangles formed by~$p_i$ and any two neighbors adjacent
to~$p_i$.  %% From~\cite{HGB:97}, 
The circumcenter of the triangle determined by $p_{i}$, $p_{j}$, and $p_k$ is
\begin{multline}\label{eq:vertices-neighbors}
  \frac{1}{4 M} \Big( \|\alpha_{kj}\|^2( \alpha_{ji} \cdot \alpha_{ik}) p_i
  + \|\alpha_{ik} \|^2(\alpha_{kj} \cdot \alpha_{ji}) p_j \\
  + \|\alpha_{ji}\|^2( \alpha_{ik} \cdot \alpha_{kj}) p_k \Big) \, ,
\end{multline}
where $M$ is the area of the triangle, and $\alpha_{ls}
=p_l-p_s$.

Equation~\eqref{eq:closed-expressions} for a polygon's
centroid and equation~\eqref{eq:vertices-neighbors} for the
Voronoi cell's vertexes lead to a closed-form
\emph{algebraic} expression for the control law in
equation~\eqref{eq:continuous-lloyd} as a function of the
neighboring vehicles' location.

\subsection{Numerical simulations} 
\label{se:numerical-simulations}
\begin{figure*}[htbp]
  \centering%  
  \includegraphics[width=.327\linewidth]{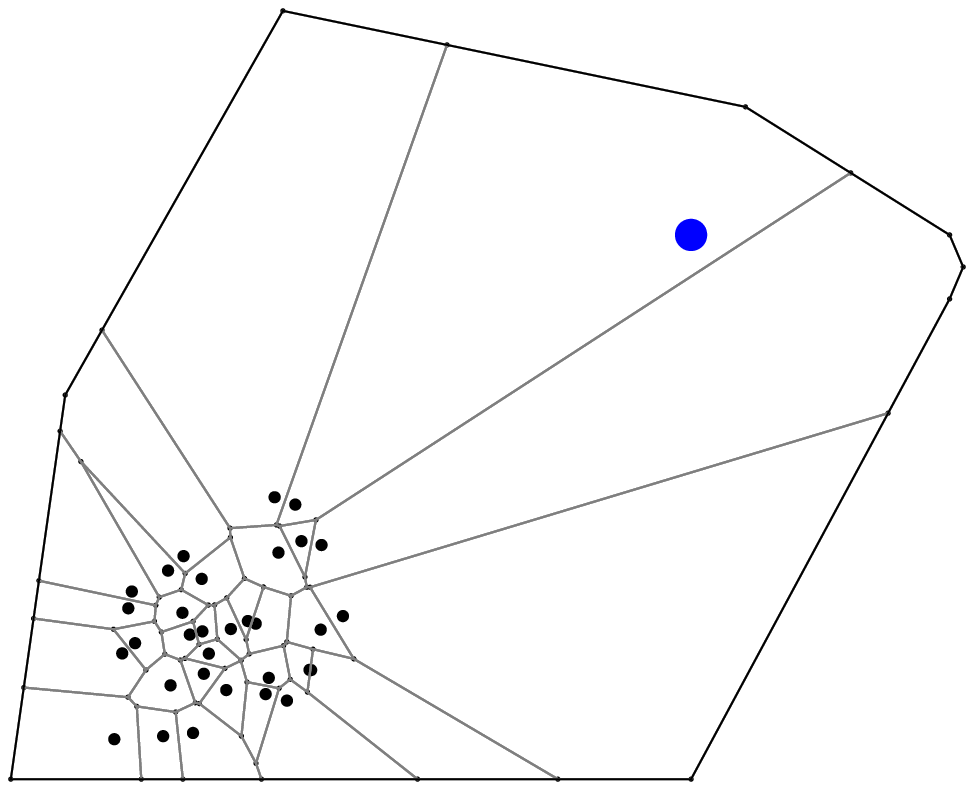} 
  \includegraphics[width=.327\linewidth]{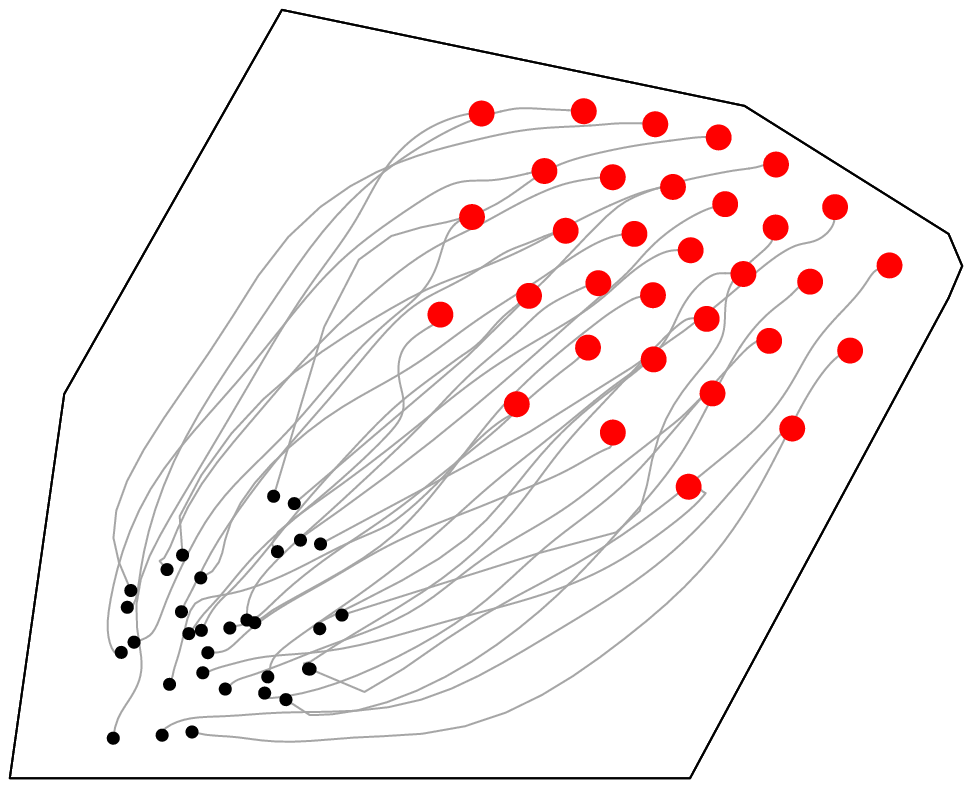}
  \includegraphics[width=.327\linewidth]{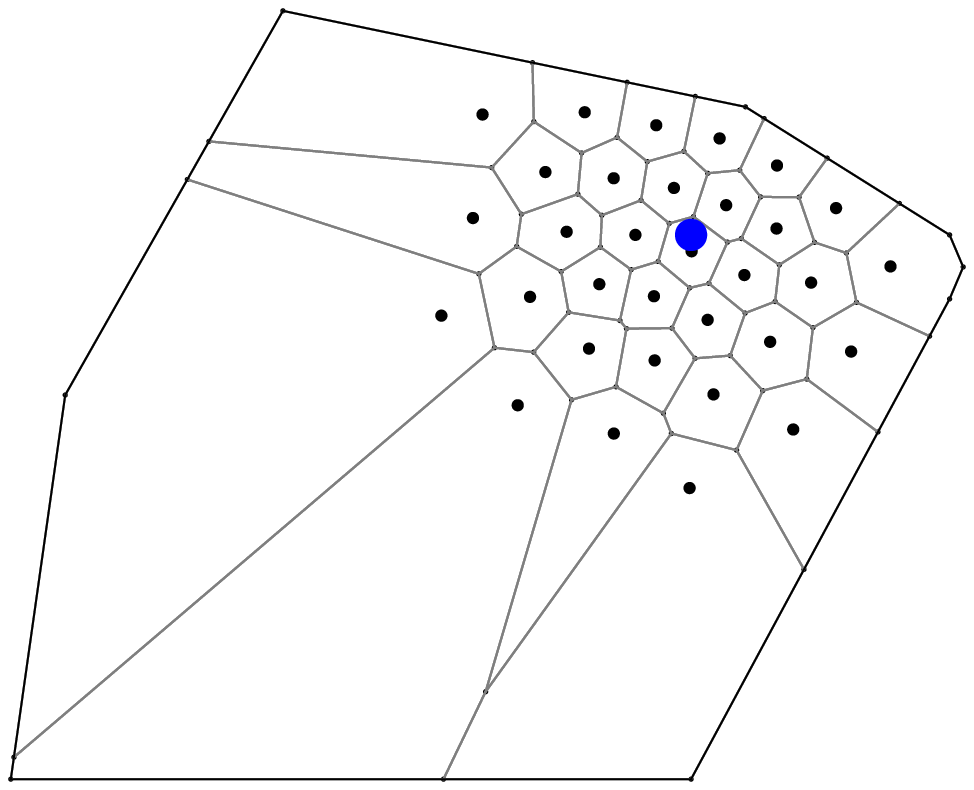}
\caption{Lloyd continuous-time algorithm for $32$ agents 
  on a convex polygonal environment, with Gaussian density $\phi=\exp(5.(-x^2
  - y^2))$ centered about the gray point in the figure.  The control gain
  in~\eqref{eq:continuous-lloyd} is $\subscr{k}{prop} = 1$ for all the
  vehicles. The left (respectively, right) figure illustrates the initial
  (respectively, final) locations and Voronoi partition. The central figure
  illustrates the gradient descent flow. }
  \label{fig:coverage-1}  
\end{figure*}

To illustrate the performance of the continuous-time Lloyd
algorithm, we include some simulation results. The algorithm
is implemented in \texttt{Mathematica} as a single
centralized program.  For the $\real^2$ setting, the code
computes the bounded Voronoi diagram using the
\texttt{Mathematica} package \texttt{ComputationalGeometry},
and computes mass, centroid, and polar moment of inertia of
polygons via the numerical integration
routine~\texttt{NIntegrate}.  Careful attention was paid to
numerical accuracy issues in the computation of the Voronoi
diagram and in the integration.  We illustrate the
performance of the closed-loop system in
Figure~\ref{fig:coverage-1}.

\section{Asynchronous distributed implementations}
\label{sec:distributed-asynchronous} 

In this section we show how the Lloyd gradient algorithm can be implemented
in an asynchronous distributed fashion. In Section~\ref{subsec:modeling} we
describe our model for a distributed asynchronous network of robotic agents.
Next, we provide two distributed algorithms for the local computation and
maintenance of the Voronoi cells. Finally, in
Section~\ref{subsec:asynchronous-Lloyd} we propose two distributed
asynchronous implementations of Lloyd algorithm: the first one is based on
the gradient optimization algorithms as described in~\cite{JNT-DPB-MA:86} and
the second one relies on the special structure of the coverage problem.

\subsection{Modeling an asynchronous distributed network of mobile robotic
  agents}\label{subsec:modeling}

We start by modeling a robotic agent that performs sensing, communication,
computation, and control actions.  We are interested in the behavior of the
asynchronous network resulting from the interaction of finitely many robotic
agents.  A theoretical framework to formalize the following concepts is that
developed in the theory of distributed algorithms; see~\cite{NAL:97}.

Let us here introduce the notion of \emph{robotic agent with computation,
  communication, and control capabilities} as the $i$th element of a network.
The $i$th agent has a processor with the ability of allocating continuous and
discrete states and performing operations on them.  Each vehicle has access
to its unique identifier $i$.  The $i$th agent occupies a location $p_i \in Q
\subset \real^N$\ and it is capable of moving in space, at any time
$t\in\real_+$ for any period of time $\delta t\in\real_+$, according to a
first order dynamics of the form:
\begin{equation*}
  \dot{p}_i(s) = u_i, \qquad \|u_i\| \leq 1 \, , \quad \forall s\in
  [t,t+\delta t]. 
\end{equation*}
The processor has access to the agent's location $p_i$ and determines the
control pair $(\delta t, u_i)$.  The processor of the $i$th agent has access
to a local clock $t_i\in\real_+\cup\{0\}$, and a \emph{scheduling sequence},
i.e., an increasing sequence of times
$\setdef{T_{i,k}\in\real_+\cup\{0\}}{k\in\natural\cup\{0\}}$ such that
$T_{i,0}=0$ and $\subscr{t}{$i$,min}<T_{i,k+1}-T_{i,k}<\subscr{t}{$i$,max}$.
The processor of the $i$th agent is capable of transmitting information to
any other agent within a closed disk of radius $R_i\in\real_+$.  We assume
the communication radius $R_i$ to be a quantity controllable by the $i$th
processor and the corresponding communication bandwidth to be limited.

We shall alternatively consider networks of \emph{robotic agents with
  computation, sensing, and control capabilities}. In this case, the
processor of the $i$th agent has the same computation and control
capabilities as before. Furthermore, we assume the processor can detect any
other agent within a closed disk of radius $R_i\in\real_+$.  We assume the
sensing radius $R_i$ to be a quantity controllable by the processor.

\subsection{Voronoi cell computation and maintenance} 
% Note that the partial derivative of $\HH_{\VV}$ with respect to the position
% of the sensor $i$ depends uniquely on the Voronoi polytope $V_i$, which in
% turn depends on the position of its Voronoi neighbors.  We therefore say that
% the computation of the derivative is completely decentralized \emph{in the
%   sense of Voronoi}.

A key requirement of the Lloyd algorithms presented in
Section~\ref{sec:coverage-control} is that each agent must be able to compute
its own Voronoi cell.  To do so, each agent needs to know the relative
location (distance and bearing) of each Voronoi neighbor.  The ability of
locating neighbors plays a central role in numerous algorithms for
localization, media access, routing, and power control in ad-hoc wireless
communication networks; e.g.,
see~\cite{JG-LJG-JH-LZ-AZ:01,XYL-PJW:01a,SM-SS-VK-MP:01,MC-CH:02a} and
references therein.  Therefore, any motion control scheme might be able to
obtain this information from the underlying communication layer.  In what
follows, we set out to provide a distributed asynchronous algorithm for the
local computation and maintenance of Voronoi cells.  The algorithm is related
to the synchronous scheme in~\cite{MC-CH:02a} and is based on basic
properties of Voronoi diagrams.

We present the algorithm for a robotic agent with sensing capabilities (as
well as computation and control).  The processor of the $i$th agent allocates
the information it has on the position of the other agents in the state
variable $P^i$. The objective is to determine the smallest distance $R_i$ for
vehicle $i$ which provides sufficient information to compute the Voronoi cell
$V_i$. We start by noting that $V_i$ is a subset of the convex set
\begin{align} \label{eq:superset-voronoi-polygon}
  W(p_i,R_i) &= \overline{B}(p_i,R_i) \cap 
  \big(\cap_{j: \|p_i-p_j\| \leq R_i} S_{ij}\big),
\end{align}
where $\overline{B}(p_i,R_i) = \setdef{q\in Q}{\| q-p_i \| \le R_i}$ and the
half planes $S_{ij}$ are
\begin{gather*}
  \setdef{q\in\real^N}{2 q\cdot(p_i-p_j)\geq (p_i+p_j)\cdot(p_i-p_j)}.
\end{gather*}
Provided $R_i$ is twice as large as the maximum distance between $p_i$ and
the vertexes of $W(p_i,R_i)$, all Voronoi neighbors of $p_i$ are within
distance $R_i$ from $p_i$ and the equality $V_i=W(p_i,R_i)$ holds.  The
minimum adequate sensing radius is therefore $R_{i,\min}=2\max_{q\in
  W(p_i,R_{i,\min})}\|p_i-q\|$.  We are now ready to state the following
algorithm.

\smallskip
\noindent{\framebox[\linewidth]{\noindent\parbox{\linewidth-2\fboxsep}{%
      \noindent\begin{tabular}{ll}
        \textbf{Name:}      & \algoASR\\
        \textbf{Goal:}      &
        \parbox[t]{.67\linewidth}{distributed Voronoi cell}
        \\
        \textbf{Requires:}      &
        \parbox[t]{.67\linewidth}{sensor with  radius $R_i$} \\[0.2ex] \hline
      \end{tabular}\\[.5ex]
      
      Local agent $i$ performs: \null\hfill\null
      \begin{algorithmic}[1]
        \STATE initialize $R_i$, detect vehicles $p_j$ within radius $R_i$
        \STATE update $P^i(t_i)$, compute $W(p_i(t_i),R_i)$
        \WHILE{  $R_{i} <  2\max_{q\in
            W(p_i(t_i),R_{i})}\|p_i(t_i)-q\|$ }
        \STATE $R_i := 2 \max_{q\in W(p_i(t_i),R_{i})}\|p_i(t_i)-q\|$
        \STATE detect vehicles $p_j$ within radius $R_i$
        \STATE update $P^i(t_i)$
        \STATE compute $W(p_i(t_i),R_i)$
        \ENDWHILE
      \STATE set $R_i := 2 \max_{q\in
        W(p_i(t_i),R_{i})}\|p_i(t_i)-q\|$ 
      \STATE set $V_i := W_i(p_i(t_i),R_i)$
      \end{algorithmic} }}}
\medskip

A similar algorithm can be designed for a robotic agent with communication
capabilities. The specifications go as in the previous algorithm, except for
the fact that steps {\small \texttt{2:}} and {\small \texttt{7:}} are
substituted by
\begin{flushleft}
\quad send $\big(\text{``request to reply''},p_i(t_i)\big)$ within radius $R_i$
     
\quad receive $\big(\text{``response''},p_j\big)$ from all agents within radius
    $R_i$ 
\end{flushleft}
Further, we have to require each agent to perform the following event-driven
task: if the $i$th agent receives at any time $t_i$ a ``request to reply''
message from the $j$th agent located at position $p_j$, it executes
\begin{center}
  send $\big(\text{``response''},p_i(t_i)\big)$  within radius $\|p_i(t)-p_j\|$
\end{center}
We call this algorithm \algoACR.

Next, we present an algorithm whose objective is to maintain the information
about the Voronoi cell of the $i$th agent, and detect the presence of certain
events.  We consider only robotic agents with sensing capabilities.  We call
an agent active if it is moving and we assume that the $i$th agent can
determine if any agent within radius $R_i$ is active or not.  Two events are
of interest: (i) a Voronoi neighbor of the $i$th agent becomes active and
(ii) a new active agent becomes a Voronoi neighbor of the $i$th agent. In
both cases, we require a trigger message ``request recomputation'' to an
appropriate control algorithm that we shall present in the next section.
Before presenting the algorithm, let us introduce the map
$\operatorname{weight}$ that assigns to the state vector
$P^i\in\real^{N\times{n}}$ a tuple $(w_1,\dots,w_n)\in\natural^n$ according
to
\[
w_j = 
\left\{
\begin{array}{cl}
3 & \text{if $j \in \NN (i)$ and $j$ is active}\\
1 & \text{if $j \in \NN (i)$ and $j$ is not active}\\
0 & \text{if $j \not \in \NN (i)$} \,.
\end{array}
\right.
\]
The algorithm is designed to run for times $t_i\in [t_0,t_0+\delta t]$.

\smallskip
\noindent{\framebox[\linewidth]{\noindent\parbox{\linewidth-2\fboxsep}{%
      \noindent\begin{tabular}{ll}
        \textbf{Name:}      & \algoM\\
        \textbf{Goal:}      &
        \parbox[t]{.67\linewidth}{Cell maintenance \& event
        detection}
        \\
        \textbf{Requires:}      &
        \parbox[t]{.67\linewidth}{
          (i) sensor with  radius $R_i$\\
          (ii) positive reals $t_0$, $\delta t$\\ 
          (iii) \algoASR
          } \\[8ex] \hline
      \end{tabular}\\[.5ex]
      
      Local agent $i$ performs for $t_i\in [t_0,t_0+\delta t]$:\null
      \begin{algorithmic}[1]
        \STATE initialize $P^i(t_0)$ and $V_i(t_0)$, set $w=
        \text{weight}(P^i(t_0))$
        \WHILE{$t_i \le t_0 + \delta t$}
        \STATE run \algoASR\
        \IF{$\text{weight}_j(P^i(t_i)) \geq w_j +2$}
        \STATE send (``request recomputation'')
        \STATE set $w=\text{weight}(P^i(t_i))$
        \ENDIF
        \ENDWHILE
      \end{algorithmic} }}}
\medskip

\subsection{Asynchronous distributed implementations of coverage control}
\label{subsec:asynchronous-Lloyd}

Let us now present two versions of Lloyd algorithm for the solution of the
optimization problem~\eqref{eq:HH} that can be implemented by an asynchronous
distributed network of robotic agents.  For simplicity, we assume that at
time $0$ all clocks are synchronized (although they later can run at
different speeds) and that each agent knows at $0$ the exact location of
every other agent.  The first algorithm is designed for robotic agents with
communication capabilities, and requires the \algoACR\, (while it does not
require the \algoM).

\smallskip

\noindent\framebox[\linewidth]{\noindent\parbox{\linewidth-2\fboxsep}{%
    \noindent\begin{tabular}{ll}
      \textbf{Name:} & \algoCone \\
      \textbf{Goal:} &
      \parbox[t]{.67\linewidth}{distributed optimal agent location} \\
%      \textbf{Assumes:}      & 
%      \parbox[t]{.67\linewidth}{${p}_i(t+1)=p_i(t)+u_i$, $\|u_i\|\leq 1$} \\
      \textbf{Requires:} & \parbox[t]{.67\linewidth}{
        (i) Voronoi cell computation\\
        (ii) centroid and mass computation \\
        (iii) positive real $\delta_0$ \\
        (iv) \algoACR}\\[11ex] \hline
    \end{tabular}\\[.5ex]

For  $i\in\{1,\dots,n\}$, $i$th agent performs at $t_i=T_{i,0}=0$: \\[-2ex]

    \begin{algorithmic}[1]
    \item $P^i(T_{i,0}) := (p^i_1(T_{i,0}),\dots,p^i_n(T_{i,0}))$
    \item compute Voronoi region $V_i(T_{i,0})$
    \item set $V_i = V_i(T_{i,0})$ and $R_i = 2 \max_{q\in V_i} \|p_i-q\|$
    \end{algorithmic}

    \medskip
    
    For $i\in\{1,\dots,n\}$, the $i$th agent performs at time $t_i=T_{i,k}$
    either one of the following threads or both.  For some $B_i\in\natural$,
    we require that after $B_i$ steps of the scheduling sequence, each of the
    threads has been executed at least once.  \\[-2ex]

    \medskip

[Information thread]\\[-3ex]
\begin{algorithmic}[1]
  \STATE run \algoACR
\end{algorithmic}
\smallskip
[Control thread]\\[-3ex]
    \begin{algorithmic}[1]
      \STATE compute centroid $C_{V_i}$ and mass $M_{V_i}$ of $V_i$
      \STATE apply control pair $\big(\delta_0, \; M_{V_i}(C_{V_i} -
      p_i(T_{i,k})) \big)$ 
    \end{algorithmic}
    }}

\medskip

As a consequence of  Theorem~3.1 and Corollary~3.1 in~\cite{JNT-DPB-MA:86}, we
have the following result.
\begin{proposition}
  Let $P_0 \in Q^n$ denote the initial sensors location.  Let $\{T_k\}$ be
  the sequence in increased order of all the scheduling sequences of the
  agents of the network.  Assume $\inf_k \{ T_k - T_{k-1}\}>0$.  Then, there
  exists a sufficiently small $\delta_* >0$ such that if $0 < \delta_0 \le
  \delta_*$, the \algoCone\, converges to the set of critical points of
  $\HH_\VV$, that is, the set of centroidal Voronoi configurations.
\end{proposition}

Next, we focus on distributed asynchronous implementations of Lloyd algorithm
that take advantage of the special structure of the coverage problem.  The
following algorithm is designed for robotic agents with sensing capabilities,
it requires the Monitoring and the Adjust sensing radius algorithms. Two
advantages of this algorithm over the previous one are that there is no need
for each agent to exactly go toward the centroid of its Voronoi cell nor to
take a small step at each stage.

\smallskip
\noindent\framebox[\linewidth]{\noindent\parbox{\linewidth-2\fboxsep}{%
    \noindent\begin{tabular}{ll}
      \textbf{Name:}      & \algoCtwo\\
      \textbf{Goal:}      & 
      \parbox[t]{.67\linewidth}{distributed optimal agent location}  
      \\
      \textbf{Requires:}      & 
      \parbox[t]{.67\linewidth}{
        (i) Voronoi cell computation\\
        (ii) centroid computation\\
        (iii) \algoM}\\[4ex] \hline
    \end{tabular}\\[.5ex]

For  $i\in\{1,\dots,n\}$, $i$th agent performs at $t_i=T_{i,0}=0$: \\[-2ex]

    \begin{algorithmic}[1]
    \item $P^i(T_{i,0}) := (p^i_1(T_{i,0}),\dots,p^i_n(T_{i,0}))$
    \item compute Voronoi region $V_i(T_{i,0})$
    \item set $V_i = V_i(T_{i,0})$ and $R_i = 2 \max_{q\in V_i} \|p_i-q\|$
    \end{algorithmic}
    
For $i\in\{1,\dots,n\}$, $i$th agent performs at $t_i=T_{i,k}$: \\[-2ex]
    \begin{algorithmic}[1]
      \STATE choose $0<\delta t_i < t_{i,min}$
      \STATE set $s=T_{i,k}$, compute  centroid $C_{V_i}(s)$
      \STATE choose  $u_i$,  with  $u_i \cdot (C_{V_i}-p_i(s)) \ge  0$, with
      strict inequality if  $p_i (s)\not = C_{V_i}$
      \WHILE{$ t_i \le T_{i,k} + \delta t_i$}
      \STATE run \algoM\, for $(T_{i,k},\delta t_i)$
      \WHILE{no warning}
      \STATE $\dot{p_i} = u_i$
      \ENDWHILE
      \STATE set $s=t_i$, compute centroid $C_{V_i(s)}$
      \STATE choose  $u_i$,  with  $u_i \cdot (C_{V_i}-p_i(s)) \ge  0$, with
      strict inequality if  $p_i (s)\not = C_{V_i}$      
      \ENDWHILE
    \end{algorithmic} }}

\medskip

% Anytime a resource is shared between users, a notion of fairness arise
% \cite{SHUL-DEL:99}. It might be easy to add logic to the control law
% whereby vehicles with higher benefit negotiate relative weights of
% the Gaussian functions they are currently covering.

\begin{remark}
  The control law $u_i$ in step {\small \texttt{7:}} can be defined via a
  saturation function. For instance, $\SR: \real^N \rightarrow \real^N$
  \[
  \SR(x) = \left\{
    \begin{array}{cc}
      x & \hbox{if $\|x\| \le 1$}\\
      x/\|x\| & \hbox{if $\|x\| \ge 1$}
    \end{array}
  \right.
  \]
  Then set $u_i = \SR (C_{V_i}-p_i)$.
\end{remark}

Resorting to the discussion in Section~\ref{se:discrete-Lloyd} on the
convergence of the discrete Lloyd algorithms, one can prove that the Coverage
behavior algorithm II verifies properties (a) and (b). As a consequence
of Proposition~\ref{prop:discrete-Lloyd}, we then have the following result.

\begin{proposition}
  Let $P_0 \in Q^n$ denote the initial sensors location.  The \algoCtwo\,
  converges to the set of critical points of $\HH_\VV$, that is, the set of
  centroidal Voronoi configurations.
\end{proposition}

\section{Extensions and applications}\label{sec:extensions-applications}

In this section we investigate various extensions and applications of the
algorithms proposed in the previous sections. We extend the treatment to
vehicles with passive dynamics and we also consider discrete-time
implementations of the algorithms for vehicles endowed with a local motion
planner. Finally, we describe interesting ways of designing density functions
to solve problems apparently unrelated to coverage.

\subsection{Variations on vehicle dynamics}\label{sec:vehicle-dynamics}

Here, we consider vehicles systems described by more general linear and
nonlinear dynamical models.
 
\emph{Coordination of vehicles with passive dynamics}.  We start by
considering the extension of the control design to nonlinear control systems
whose dynamics is passive.  Relevant examples include networks of vehicles
and robots with general Lagrangian dynamics, as well as spatially invariant
passive linear systems.  Specifically, assume that for each
$i\in\{1,\dots,n\}$, the $i$th vehicle state includes the spatial variable
$p_i$, and that the vehicle's dynamics is passive with input~$u_i$,
output~$\dot{p}_i$ and storage function $S_i:Q \rightarrow \real_+$.
Furthermore, assume that the input preserving the zero dynamics manifold $\{
\dot{p}_i = 0\}$ is $u_i = 0$.

%\todo{guarantees on not leaving Q?}

For such systems, we devise a proportional derivative (PD) control via,
\begin{gather}\label{eq:second-order}
  u_i= - \subscr{k}{prop} M_{V_i} (p_i - C_{V_i}) - \subscr{k}{deriv} \dot{p}_i,
\end{gather}
where $\subscr{k}{prop}$ and $\subscr{k}{deriv}$ are scalar positive gains.
The closed-loop system induced by this control law can be analyzed with the
Lyapunov function
\begin{equation*}
  {\mathcal{E}} = \frac{1}{2} \subscr{k}{prop}\HH_{\VV}+\sum_{i=1}^n S_i,
\end{equation*}
yielding the following result.

\begin{proposition} \label{le:second-order}
  For passive systems, the control law~\eqref{eq:second-order} achieves
  asymptotic convergence of the sensors location to the set of centroidal
  Voronoi configurations. If this set is finite, then the sensors location
  converges to a centroidal Voronoi configuration.
\end{proposition}

\begin{proof}
Consider the evolution of the function ${\mathcal{E}}$,
\begin{multline*}
  \frac{d}{dt} {\mathcal{E}} = \frac{1}{2} k_{\text{prop}} \frac{d}{dt} \HH_{\VV} +
  \sum_{i=1}^n \dot{S}_i \\
  \le \subscr{k}{prop} M_{V_i} (p_i - C_{V_i}) + \dot{p}_i u_i = -\subscr{k}{deriv}
  \sum_{i=1}^n\dot{p}_i^2 \le 0 \, .
\end{multline*}
By LaSalle's principle, the sensors location converges to the largest
invariant set contained in $\{ \dot{p}_i = 0 \}$. Given the assumption on the
zero dynamics, we conclude that $p_i = C_{V_i}$ for $i \in \{ 1, \dots,n \}$,
i.e., the largest invariant set corresponds to the set of centroidal Voronoi
configurations.  If this set is finite, LaSalle's principle also guarantees
convergence to a specific centroidal Voronoi configuration.
\end{proof}

\begin{figure}[htbp]
  \centering
  \includegraphics[width=.33\linewidth]{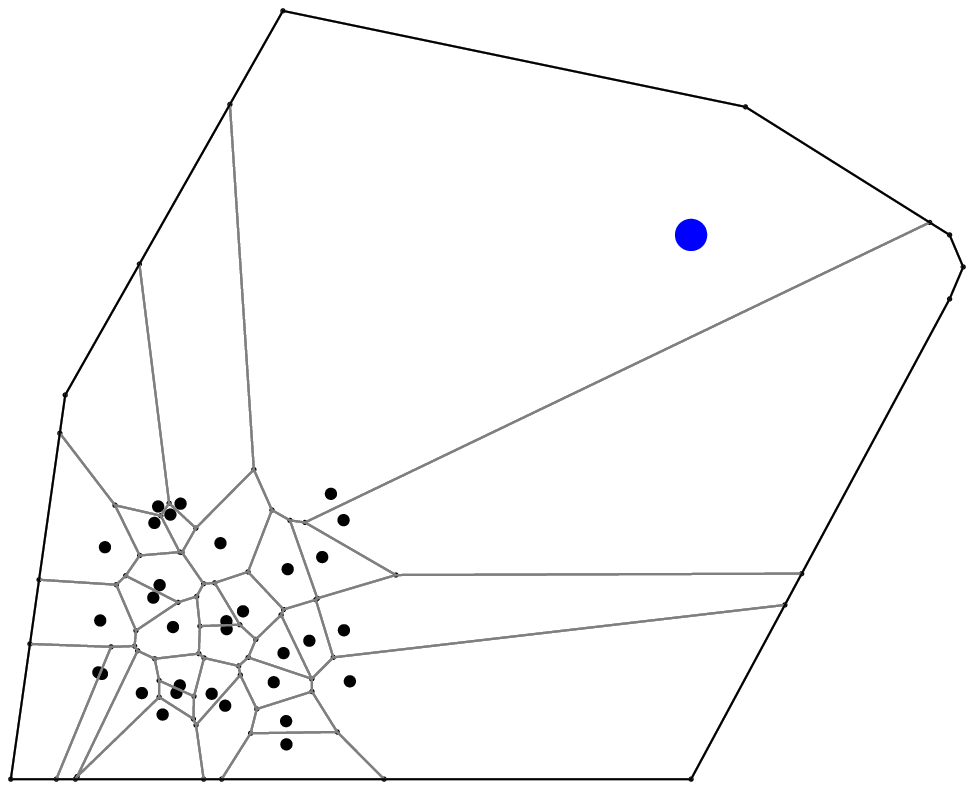}%
  \includegraphics[width=.33\linewidth]{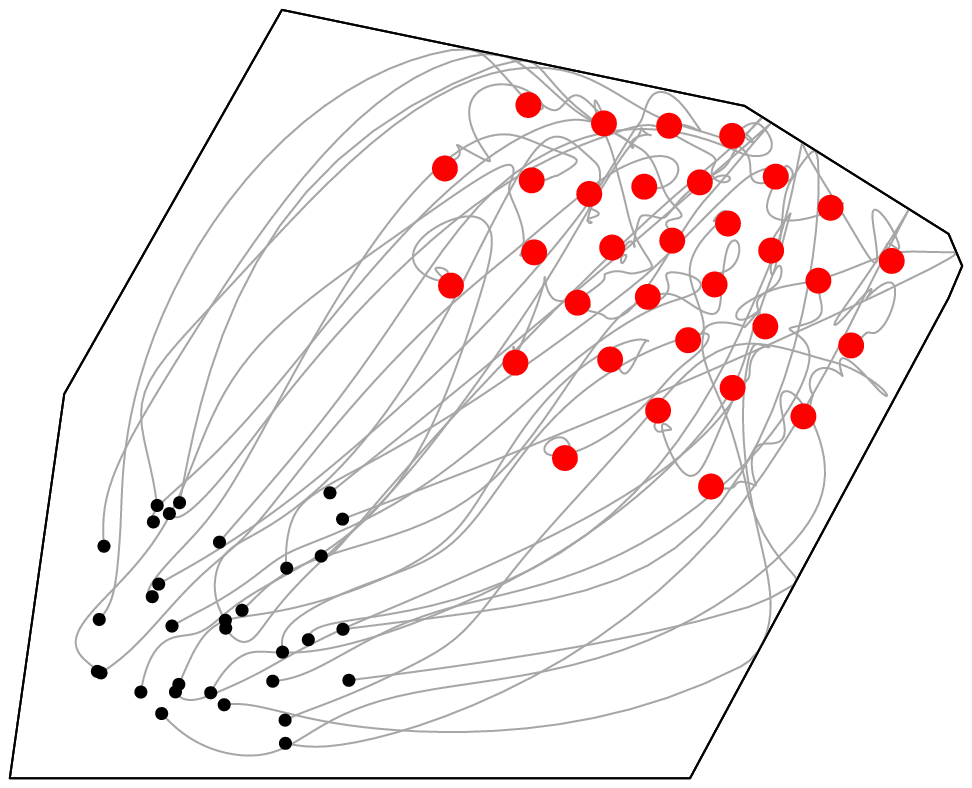}%
  \includegraphics[width=.33\linewidth]{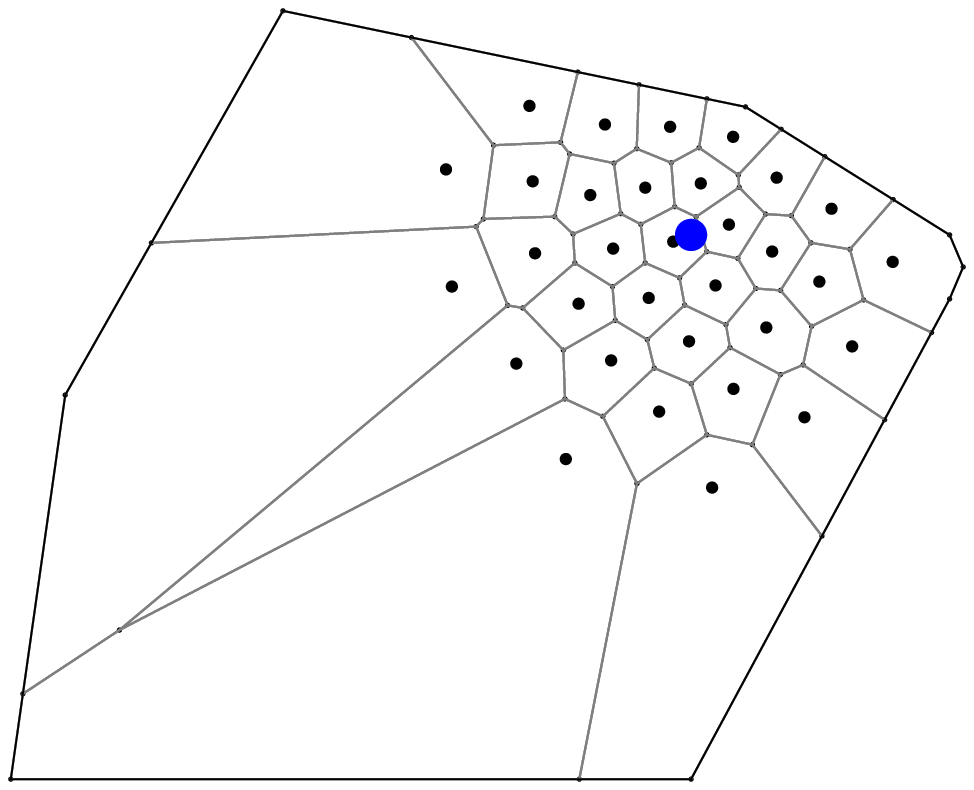}%
  \caption{Coverage control for $32$ vehicles with second order dynamics.
    The environment and Gaussian density function are as in
    Figure~\ref{fig:coverage-1}. The control gains are $\subscr{k}{prop} = 6$
    and $\subscr{k}{deriv}=1$.}
  \label{fig:second-order}
\end{figure}

In Figure~\ref{fig:second-order} we illustrate the performance of the control
law~\eqref{eq:second-order} for vehicles with second-order dynamics
$\ddot{p}_i=u_i$.

\emph{Coordination of vehicles with local controllers}.  Next, consider the
setting where each vehicle has an arbitrary dynamics and is endowed with a
local feedback and feedforward controller. The controller is capable of
strictly decreasing the distance to any specified position in $Q$ in a
specified period of time $\delta$.  

Assume the dynamics of the $i$th vehicle is described by $\dot{x}_i = f_i
(t,x_i,u)$, where $x_i \in \real^m$ denotes its state, and $\pi_i:
\real^{m_i} \rightarrow Q$ is such that $\pi_i (x_i) = p_i$. Assume also that
for any $\subscr{p}{target} \in Q$ and any $x_0 \in \real^m \setminus
\pi_i^{-1}(\subscr{p}{target})$, there exists $\ov{u}
(t,x(t),\subscr{p}{target})$ such that the solution $\ov{x}_i(t)$ of
\[
\dot{\ov{x}}_i = f_i (t,\ov{x}_i(t),\ov{u}
(t,\ov{x}_i(t),\subscr{p}{target})) \, , \quad \ov{x}_i(0) = x_0 \, ,
\]
verifies $\| \pi_i (\ov{x}_i(t_0 + \delta)) - \subscr{p}{target} \| < \| \pi_i
(\ov{x}_i(t_0)) - \subscr{p}{target} \|$.

\begin{proposition}\label{prop:local-feedback+feedforward}  
  Consider the following coordination algorithm. At time $t_k = k \delta$,
  $k\in\natural$, each vehicle computes $V_i (t_k)$ and $C_{V_i} (t_k)$;
  then, for time $t \in [t_k,t_{k+1}[$, the vehicle executes $\ov{u}
  (t,x(t),C_{V_i}(t_k))$. For this closed-loop system, the sensors location
  converges to the set of centroidal Voronoi configurations.  If this set is
  finite, then the sensors location converges to a centroidal Voronoi
  configuration.
\end{proposition}
The proof of this result readily follows from
Proposition~\ref{prop:discrete-Lloyd}, since the algorithm verifies
properties (a) and (b) of Section~\ref{se:discrete-Lloyd}.

As an example, we consider a classic model of mobile wheeled dynamics, the
\emph{unicycle model}.  Assume the $i$th vehicle has configuration
$(\theta_i,x_i,y_i)\in\operatorname{SE}(2)$ evolving according to
\begin{align*}
  \dot{\theta}_i = \omega_i \, , \quad 
  \dot{x}_i      = v_i \cos\theta_i \, , \quad
  \dot{y}_i      = v_i \sin\theta_i \, ,
\end{align*}
where $(\omega_i,v_i)$ are the control inputs for vehicle $i$. Note that the
definition of $(\theta_i,v_i)$ is unique up to the discrete action
$(\theta_i,v_i)\mapsto(\theta_i+\pi,-v_i)$. Given a target point
$\subscr{p}{target}$, we use this symmetry to require the equality
$(\cos\theta_i,\ \sin\theta_i)\cdot (p_i- \subscr{p}{target})\leq 0$ for all
time $t$.  Should the equality be violated at some time $t=t_0$, we shall
redefine $\theta_i(t_0^+)=\theta_i(t_0^-)+\pi$ and $v_i$ as $-v_i$ from time
$t=t_0$ onwards.

\begin{figure*}[htbp]
  \centering
  \includegraphics[width=.33\linewidth]{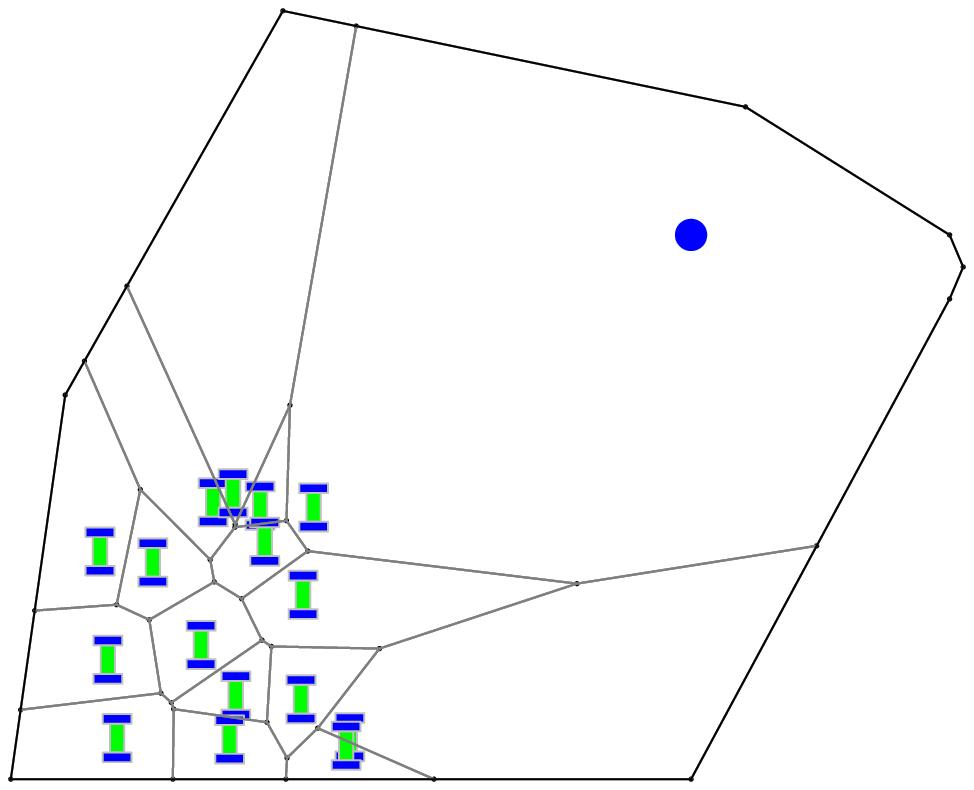}%
  \includegraphics[width=.33\linewidth]{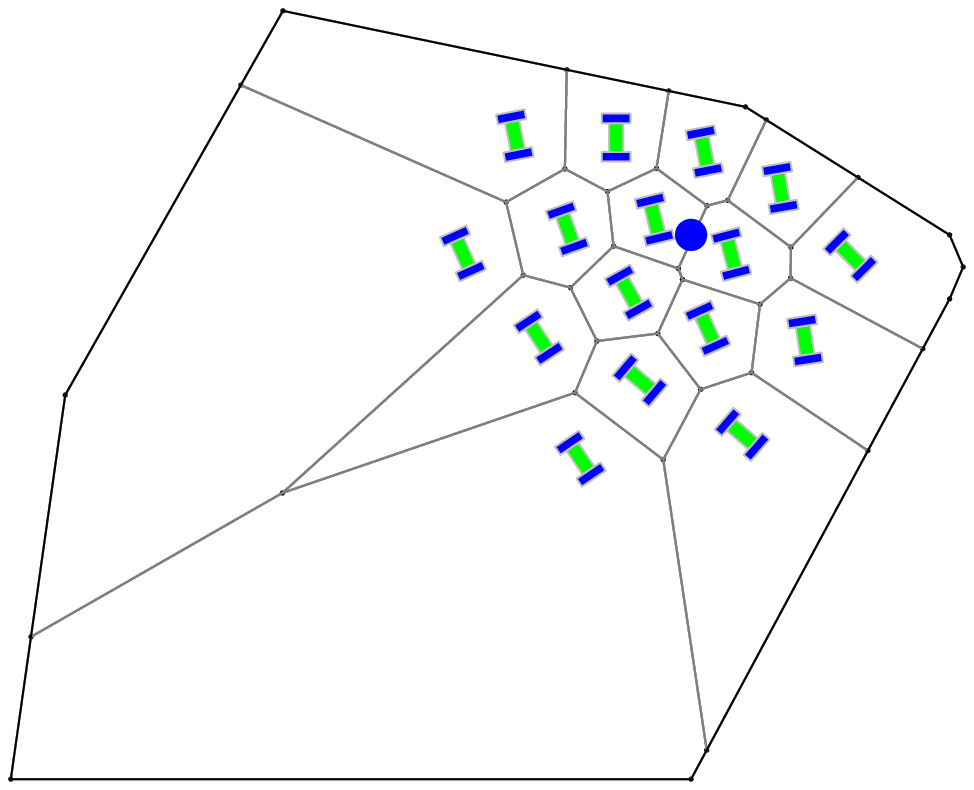}%
  \caption{Coverage control for $16$ vehicles with mobile wheeled dynamics.
    The environment and Gaussian density function are as in
    Figure~\ref{fig:coverage-1}, and $\subscr{k}{prop}=3$.}
  \label{fig:coverage-mobile}
\end{figure*}

Following the approach in~\cite{AA:99a}, consider the control law
\begin{align*}
  \omega_i &= 2 \subscr{k}{prop} \arctan \frac{(-\sin\theta_i,\ 
    \cos\theta_i)\cdot (p_i-\subscr{p}{target})}{(\cos\theta_i,\ 
    \sin\theta_i)\cdot (p_i-\subscr{p}{target})}
  \\
  v_i &= - \subscr{k}{prop} (\cos\theta_i,\ \sin\theta_i)\cdot
  (p_i-\subscr{p}{target}),
\end{align*}
where $\subscr{k}{prop}$ is a positive gain.  This feedback law differs from
the original stabilizing strategy in~\cite{AA:99a} only in the fact that no
final angular position is preferred. One can prove that $p_i = (x_i,y_i)$ is
guaranteed to monotonically approach the target position $\subscr{p}{target}$
when run over an infinite time horizon.  We illustrate the performance of the
proposed algorithm in Figure~\ref{fig:coverage-mobile}.

\subsection{Geometric patterns and formation control}
\label{sec:geometric-patterns}
Here  we suggest  the use  of decentralized  coverage algorithms  as formation
control algorithms,  and we  present various density  functions that  lead the
multi-vehicle network to predetermined  geometric patterns.  In particular, we
present simple density functions that lead to segments, ellipses, polygons, or
uniform distributions inside convex environments.

\begin{figure}[thbp]
  \centering
  \begin{tabular}{c}
    \includegraphics[width=.33\linewidth]{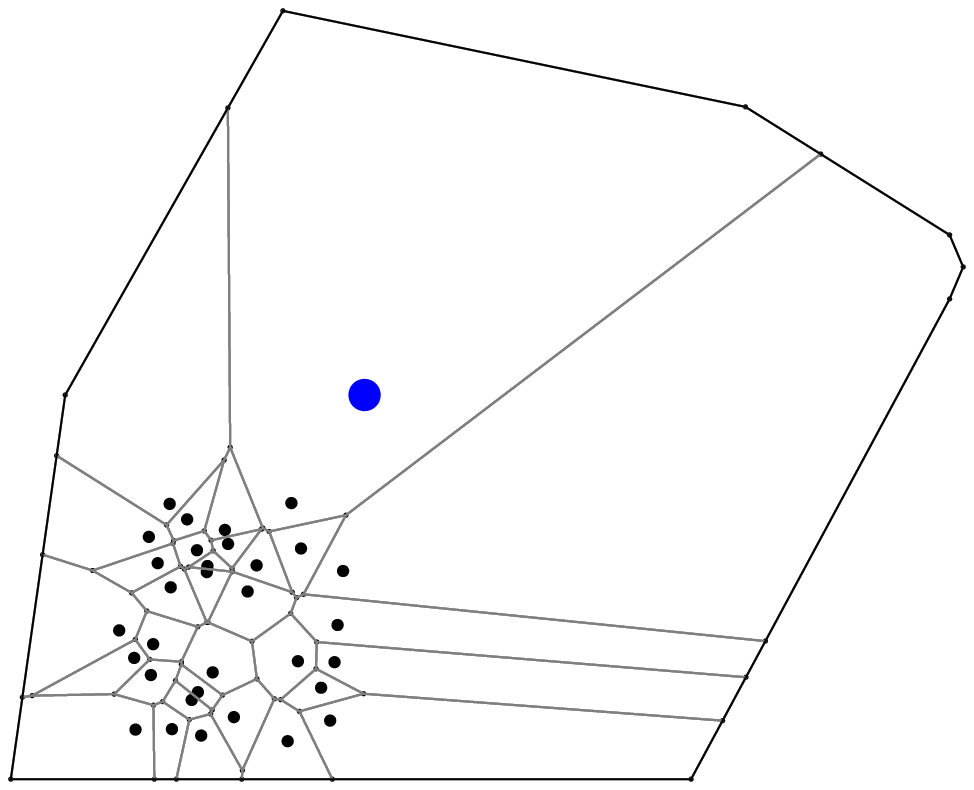}%
    \includegraphics[width=.33\linewidth]{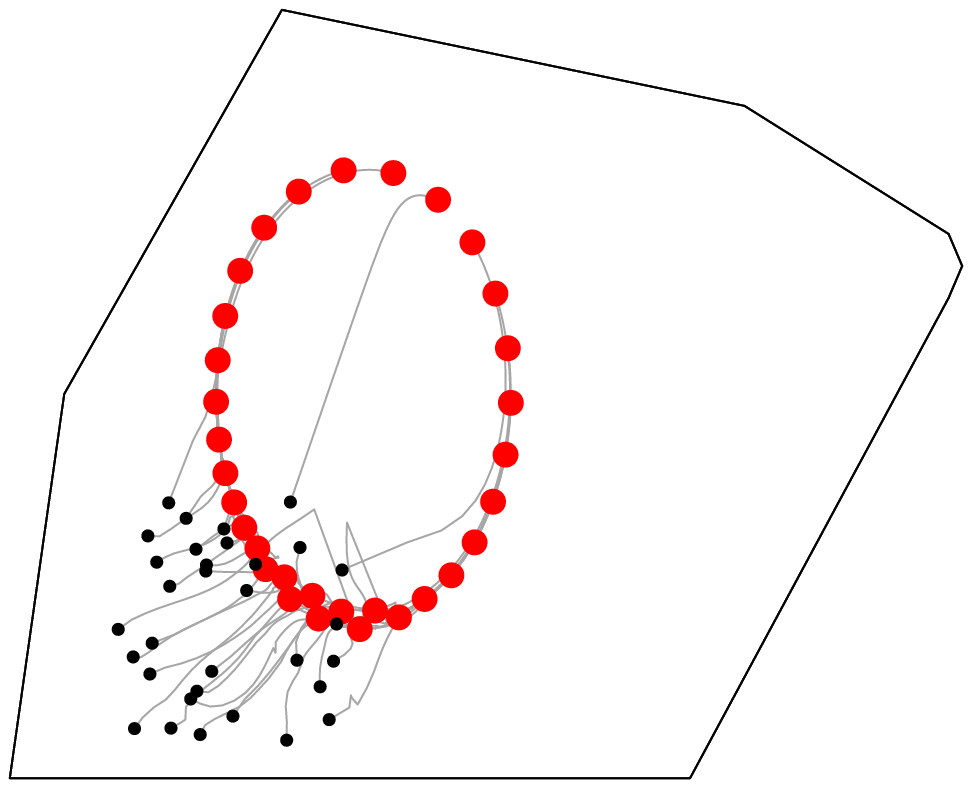}%
    \includegraphics[width=.33\linewidth]{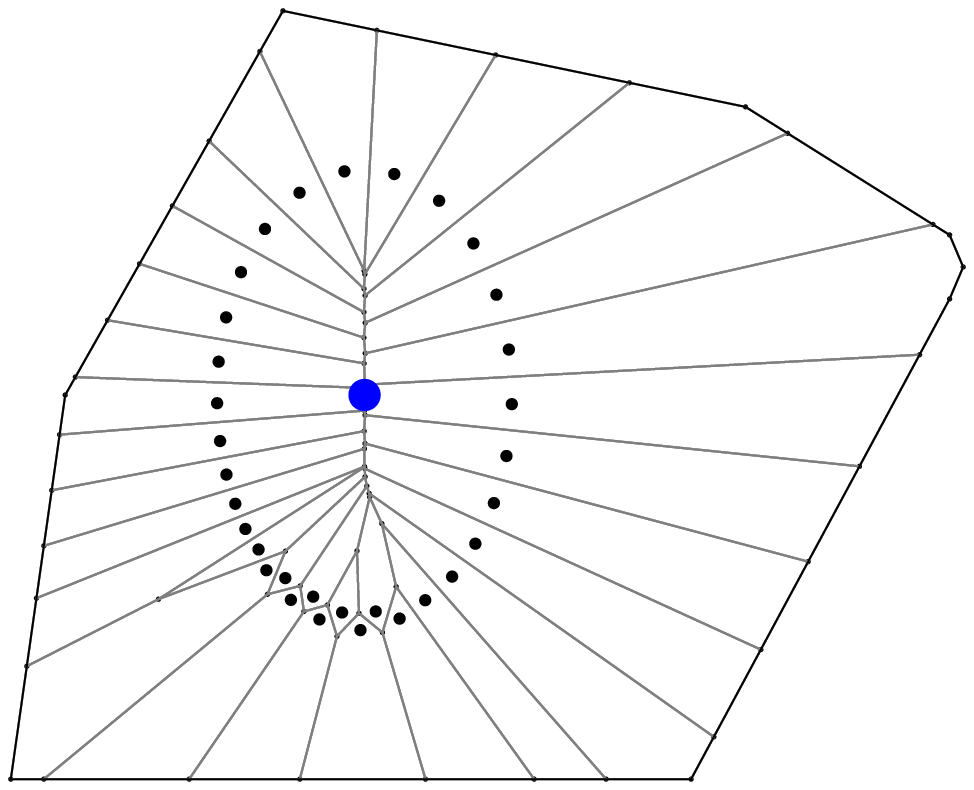}%
  \end{tabular}        
   \caption{Coverage control for $32$ vehicles with $\subscr{\phi}{ellipse}$.
     The parameter values are: $k=500$, $a=1.4$, $b=.6$, $x_c=y_c=0$,
     $r^2=.3$, and $\subscr{k}{prop}=1$.}
  \label{fig:coverage-circle}
\end{figure}

Consider a  planar environment, let $k$  be a large positive  gain, and denote
$q=(x,y)\in Q\subset \real^2$. Let $a,b,c$  be real numbers, consider the line
$ax+by+c=0$, and define the density function
\begin{align*}
  \subscr{\phi}{line}(q) &= \exp(-k(ax+by+c)^2).
\end{align*}
Similarly, let $(x_c,y_c)$  be a reference point in  $\real^2$, let $a,b,r$ be
positive  scalars, consider the  ellipse $a  (x-x_c)^2+ b  (y-y_c)^2=r^2$, and
define the density function
\begin{align*}
  \subscr{\phi}{ellipse}(q) &= \exp\big(-k (a (x-x_c)^2+ b(y-y_c)^2-r^2)^2 \big).
\end{align*}
We illustrate this density function in Figure~\ref{fig:coverage-circle}.
During the simulations, we observed that the convergence to the desired
pattern was rather slow.

\begin{figure}[thbp]
  \centering
  \includegraphics[width=.33\linewidth]{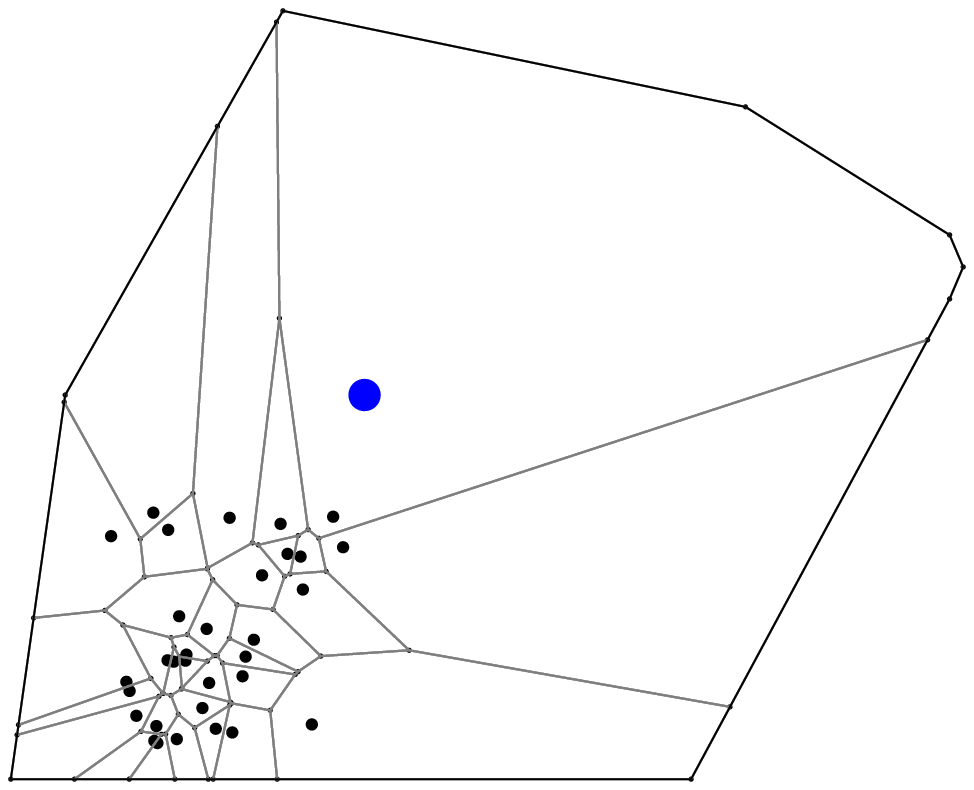}%
  \includegraphics[width=.33\linewidth]{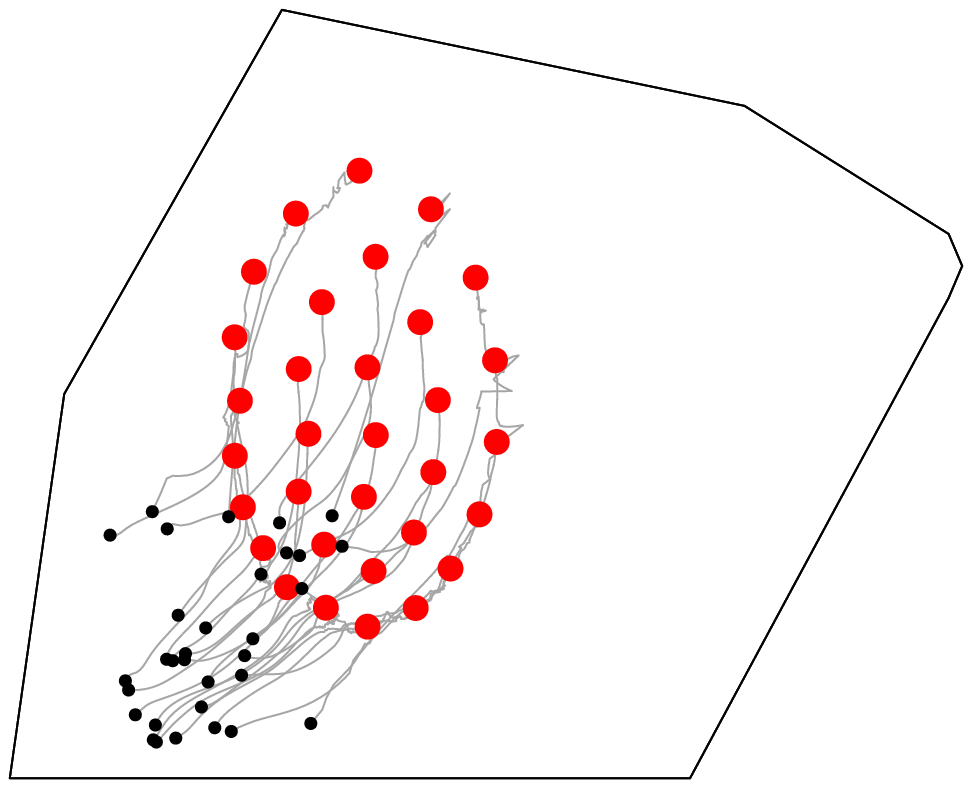}%
  \includegraphics[width=.33\linewidth]{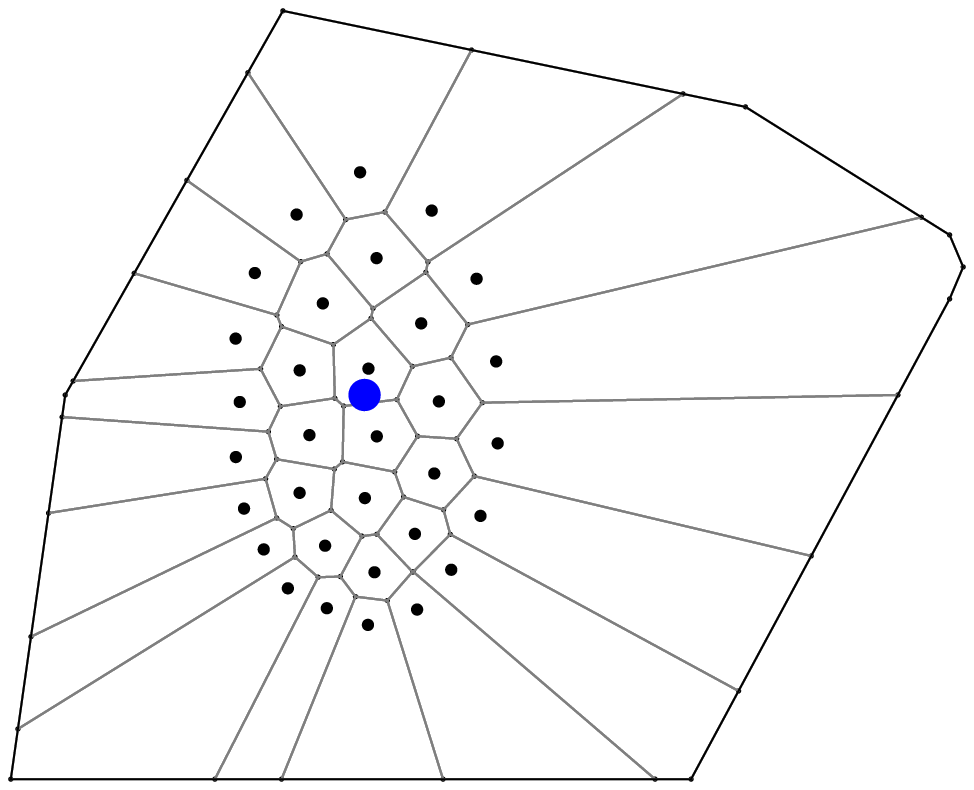}%
  \caption{Coverage control for $32$ vehicles to an ellipsoidal disk. The
    density function parameters are the same as in
    Figure~\ref{fig:coverage-circle}, and $\ell=10$, $\subscr{k}{prop}=1$.}
  \label{fig:coverage-disk}  
\end{figure}

Finally, define the smooth ramp function $\operatorname{SR}_\ell(x)= x
(\arctan(\ell x)/\pi +(1/2))$, and the density function
\begin{align*}
  \subscr{\phi}{disk}(q) =
  \exp(-k \operatorname{SR}_\ell(a(x-x_c)^2+b(y-y_c)^2-r^2)).
\end{align*}
This density function leads the multi-vehicle network to obtain a uniform
distribution inside the ellipsoidal disk $a (x-x_c)^2+b(y-y_c)^2\leq r^2$. We
illustrate this density function in Figure~\ref{fig:coverage-disk}.

It appears straightforward to generalize these types of density functions to
the setting of arbitrary curves or shapes.  The proposed algorithms are to be
contrasted with the classic approach to formation control based on rigidly
encoding the desired geometric pattern.  One disadvantage of the proposed
approach is the requirement for a careful numerical computation of Voronoi
diagrams and centroids. We refer to~\cite{KS-IS:96} for previous work on
algorithms for geometric patterns, and to~\cite{TB-RA:98,JPD-JPO-VK:01} for
formation control algorithms.

\section{Conclusions}
\label{sec:conclusions}
We have presented a novel approach to coordination algorithms for
multi-vehicle networks. The scheme can be thought of as an interaction law
between agents and as such it is implementable in a distributed asynchronous
fashion.  Numerous extensions appear worth pursuing. We plan to investigate
the setting of non-convex environments and non-isotropic sensors.  We are
currently implementing these algorithms on a network of all-terrain vehicles.
Furthermore, we plan to extend the algorithms to provide collision avoidance
guarantees and to vehicle dynamics which are not locally controllable.

% \todo{write two paragraphs}

\subsection*{Acknowledgments}
This work was supported by NSF Grant CMS-0100162, ARO Grant DAAD~190110716,
and DARPA/AFOSR MURI Award F49620-02-1-0325.

\section{Appendix}\label{sec:appendix}

In this section we collect some relevant facts on descent flows both in the
continuous and in the discrete-time settings.  We do this
following~\cite{HKK:96} and~\cite{DGL:84}, respectively.  We include
Proposition~\ref{prop:discrete-LaSalle-surprising} as we are unable to locate
it in the linear and nonlinear programming literature.

\subsection*{Continuous-time descent flows}

Consider  the differential  equation $\dot{x}  = X  (x)$, where  $X:D \subset
\real^N  \rightarrow  \real^N$  is  locally  Lipschitz and  $D$  is  an  open
connected set.  A  set $M$ is said to be  (positively) invariant with respect
to $X$ if $x(0) \in M$ implies $x(t)  \in M$, for all $t \in \real$ (resp. $t
\ge 0$).  A descent function for $X$ on $\Omega$, $V:D \rightarrow \real$, is
a continuously differentiable function such that ${\mathcal{L}}_X V \le 0$ on
$\Omega$.   We   denote  by  $E$  the   set  of  points   in  $\Omega$  where
${\mathcal{L}}_X V (x) = 0$ and by $M$ be the largest invariant set contained
in $E$.  Finally,  the distance from a point  $x$ to a set $M$  is defined as
$\dist(x,M) = \inf_{p \in M} \| x-p \|$.

\begin{lemma}[LaSalle's principle] \label{lemma:continuous-LaSalle}
Let $\Omega \subset D$ be a  compact set that it is positively invariant with
respect  to $X$.  Let  $x(0) \in  M$ and  $x_*$ be  an accumulation  point of
$x(t)$. Then $x_* \in M$  and $\dist(x(t),M) \rightarrow 0$ as $t \rightarrow
\infty$.
\end{lemma}
The following corollary is Exercise 3.22 in~\cite{HKK:96}.
\begin{corollary}\label{corollary:LaSalle-finite}
If the  set $M$ is a  finite collection of  points, then the limit  of $x(t)$
exists and equals one of them.
\end{corollary}

\subsection*{Discrete-time descent flows}

Let $X$ be  a subset of $\real^N$.  An algorithm $T$  is a continuous mapping
from $X$ to $X$. A set $C$ is said to be positively invariant with respect to
$T$ if $x_0 \in C$ implies $T(x_0) \in C$.
% For our purposes, it will be  enough to consider the case when the image of
% a point by $T$ is a single point. $T$ is  said to be \emph{closed  on
% $C$}  if  the  conditions  $x_n \rightarrow  x$,  $T(x_n) \rightarrow y$
% imply  that $y=T(x)$. For instance, if  $T$ is continuous, then it  is
% closed.
A point $x_*$ is said to be a fixed point of $T$ if
$T(x_*)=x_*$. We denote the set of fixed points of $T$ by
$\Gamma$.  A descent function for $T$ on $C$, $Z: X
\rightarrow \real_+$, is any nonnegative real-valued
continuous function satisfying $Z(T(x)) \le Z(x)$ for $x \in
C$, where the inequality is strict if $x \not \in \Gamma$.
Typically, $Z$ is the objective function to be minimized,
and $T$ reflects this goal by yielding a point that reduces
(or at least does not increase) $Z$.

\begin{lemma}[Global convergence theorem] \label{lemma:discrete-LaSalle}
Let  $C \subset X$  be a  compact set  that it  is positively  invariant with
respect to  $T$.  Let  $x_0 \in  C$ and denote  $x_m =  T (x_{m-1})$,  $m \ge
1$. Let  $x_*$ be  an accumulation point  of the  sequence $\{ x_m  \}_{m \ge
1}$.  Then $x_* \in \Gamma$, $\dist (x_m,\Gamma) \rightarrow 0$ and $Z(x_m)
\rightarrow Z(x_*)$ as $m \rightarrow \infty$.
\end{lemma}

% \begin{corollary}\label{corollary:discrete-LaSalle-singleton}
% If  the  set  $\Gamma$  consists  of  a single  point  $\bar{x}$,  then  $x_m
% \rightarrow \bar{x}$.
% \end{corollary}

% It is instructive to note that
% Corollary~\ref{corollary:discrete-LaSalle-singleton} is not valid if
% $\Gamma$ is finite and contains more that a single point. 

\begin{proposition}\label{prop:discrete-LaSalle-surprising}
If  the  set  $\Gamma$ is  a  finite  collection  of points,  then  $\{x_m\}$
converges and equals one of them.
\end{proposition}
\begin{proof}
Let $x_*$ be an accumulation point of $\{x_m\}$ and assume the whole sequence
does not converge  to it. Then, there  exists an $\eps >0$ such  that for all
$m_0$, there is a  $m' > m_0$ such that $\| x_{m'} - x_*  \| > \eps$. Let $d$
be  the minimum  of all  the distances  between the  points in  $\Gamma$. Fix
$\eps'  = \min  \{ d/2,  \eps\}$.  Since $T$  is continuous  and $\Gamma$  is
finite, there exists  $\delta > 0$ such that $\|x-z \|<  \delta$, with $z \in
\Gamma$, implies $\| T(x) - z \|  < \eps'$ (that is, for each $z \in \Gamma$,
there exists such $\delta(z)$, and we take the minimum over $\Gamma$).

Now, since $\dist  (x_m,\Gamma) \rightarrow 0$, there exists  $m_1$ such that
for all $m \ge m_1$, $\dist  (x_m,\Gamma) < \delta$. Also, we know that there
is a subsequence of $\{ x_m \}$ which converges to $x_*$, let us denote it by
$\{ x_{m_k}\}_{k \ge 1}$. For  $\delta$, there exists $m_{k_1}$ such that for
all $k \ge k_1$, we have $\| x_{m_{k}} - x_*\| < \delta$.

Let $m_0 = \max \{ m_1, m_{k_1}\}$. Take $k$ such that $m_k \ge m_0$ Then,
\begin{equation}\label{eq:inequality}
\| x_{m_{k}+1} - x_*\| = \| T(x_{m_{k}}) - x_*\| < \eps' \, .
\end{equation}
Now we are going to prove that $\| x_{m_{k}+1} - x_*\| < \delta$. If $d/2 \le
\delta$, then this claim is straightforward, since $\eps' \le d/2$. If $d/2 >
\delta$, suppose that $\| x_{m_{k}+1} -  x_*\| > \delta$. Since $m_k +1 > m_0
\ge m_1$, then $\dist (x_{m_k+1},  \Gamma) < \delta$. Therefore, there exists
$y \in \Gamma$ such that $\| x_{m_k+1}  - y \| < \delta$. Necessarily, $y \not
= x_*$. Now, by  the triangle inequality, $\| x-y \| \le  \| x-x_{m_k +1}\| +
\| x_{m_k+1} - y \|$. Then,
\begin{align*}
\| x_{m_k + 1} - x_* \| \ge \|x-y \| - \| x_{m_{k+1}} - y \| \ge d - \delta >
d/2 \, ,
\end{align*}
which contradicts~\eqref{eq:inequality}. Therefore, $\| x_{m_{k}+1} - x_*\| <
\delta$. This argument can be iterated to  prove that for all $m \ge m_0$, we
have $\|  x_m - x_*\|  < \delta$. Let  us take now $m'  > m_0$ such  that $\|
x_{m'} - x_* \| > \eps$. Since $m'-1 \ge m_0$, we have $\| x_{m'-1} - x_*\| <
\delta$, and therefore
\[
\| x_{m'} - x_* \| = \|  T(x_{m'-1}) - x_*\| < \eps' \le \eps \, ,
\]
which is a contradiction. Therefore, $\{ x_m \}$ converges to $x_*$.
\end{proof}

\end{document}